\documentclass{birkjour}
\usepackage{amsmath, amsthm, amsfonts, mathdots}
\usepackage{amssymb}
\usepackage{amsbsy}
\usepackage{mathrsfs}

\input amssym.def
\input amssym
\usepackage{enumerate}
\usepackage{hyperref}
\hypersetup{
    pdftitle={Memory Estimation of Inverse Operators},    
    pdfauthor={Anatoly Baskakov}{Ilya Krishtal},     
    pdfsubject={Memory Estimation of Inverse Operators},   
    pdfcreator={Ilya Krishtal},   
    pdfproducer={Ilya Krishtal}, 
    pdfkeywords={Banach modules}{Beurling spectrum}{Wiener's lemma}, 
     colorlinks,%
    citecolor=blue,%
    filecolor=black,%
    linkcolor=magenta,%
    urlcolor=black
}

\chardef\bslash=`\\ 

\makeatletter
\def\verbatim{\interlinepenalty\@M \@verbatim
   \leftskip\@totalleftmargin\advance\leftskip2pc
   \frenchspacing\@vobeyspaces \@xverbatim}
\makeatother
\newtheorem{thm}{Theorem}[section]
\newtheorem{cor}[thm]{Corollary}
\newtheorem{lem}[thm]{Lemma}
\newtheorem{prop}[thm]{Proposition}

\newtheorem{ass}[thm]{Assumption}

\theoremstyle{definition}
\newtheorem{defn}{Definition}[section]

\theoremstyle{remark}
\newtheorem{rem}{Remark}[section]
\newtheorem{exmp}{Example}[section]

\numberwithin{equation}{section}

\newcommand{\begeq}{\begin {equation}}
\newcommand{\eq}{\end{equation}}
\newcommand{\bs}{\begin {split}}
\newcommand{\es}{\end{split}}
\newcommand{\bp}{\begin {prop}}
\newcommand{\ep}{\end {prop}}
\newcommand{\bt}{\begin {thm}}
\newcommand{\et}{\end {thm}}
\newcommand{\bc}{\begin {cor}}
\newcommand{\ec}{\end {cor}}
\newcommand{\bl}{\begin {lem}}
\newcommand{\el}{\end {lem}}
\newcommand{\bpf}{\begin {proof}}
\newcommand{\epf}{\end {proof}}
\newcommand{\bi}{\begin {itemize}}
\newcommand{\ei}{\end {itemize}}
\newcommand{\ben}{\begin {enumerate}}
\newcommand{\een}{\end {enumerate}}
\newcommand{\brem}{\begin {rem}}
\newcommand{\erem}{\end {rem}}

\newcommand{\bd}{\begin {defn}}
\newcommand{\ed}{\end {defn}}

\newcommand{\bex}{\begin {exmp}}
\newcommand{\eex}{\end {exmp}}

\newcommand{\la}{\langle}
\newcommand{\ra}{\rangle}

\newcommand{\TTT}{{\T\kern-.44em \T}}
\newcommand{\tTTT}{\widetilde{\T\kern-.44em \T}}

\newcommand{\ZZ}{{\mathbb Z}}

\newcommand{\RR}{{\mathbb R}}
\newcommand{\CC}{{\mathbb C}}
\newcommand{\NN}{{\mathbb N}}

\newcommand{\HH}{{\mathcal H}}

\newcommand{\s}{\sigma}
\renewcommand{\l}{\lambda}
\renewcommand{\L}{\Lambda}
\renewcommand{\a}{\alpha}

\newcommand{\M}{\mathcal{M}}

\newcommand{\g}{\gamma}

\newcommand{\T}{\mathcal{T}}

\begin{document}

%
%
%
%
%
%
%
%
%


\title[General Dirac Operators as Generators of  Operator Groups]
{
General Dirac Operators as \\ Generators of   Operator Groups} 

\author{Anatoly G. Baskakov}
\address{Department of Applied Mathematics and Mechanics \\ Voronezh State University \\ Voronezh 394693 \\ Russia}
\email{anatbaskakov@yandex.ru}
\thanks{The first author is supported in part by the Ministry of Education and Science of the Russian Federation in the frameworks of the project part of the state work quota (Project No 1.3464.2017/4.6).
The second author is supported in part by NSF grant DMS-1322127. The third author is supported in part by RFBR grant 16-01-00197.}



\author{Ilya A. Krishtal}
\address{Department of Mathematical Sciences\\ Northern Illinois University\\ DeKalb, IL 60115 \\ USA}
\email{krishtal@math.niu.edu}
\author{Natalia B. Uskova}
\address{Department of Higher Mathematics and Mathematical Physical Modeling \\ Voronezh State Technical University \\ Voronezh 394026 \\ Russia}
\email{nat-uskova@mail.ru}


\date{\today }




\begin{abstract}
We use the method of similar operators to study a general Dirac operator $L$ and its spectral properties. We find a similar operator to the Dirac operator that is an orthogonal direct sum of simpler operators. The result is used to describe an operator group generated by the  operator $iL$ and study its properties such as the asymptotics of the spectrum.
\end{abstract}
\maketitle

\section{Introduction}

Dirac operators arise in various problems modeling physical phenomena such as electromagnetic fields. The operators are a centerpiece of classical boundary value problems in mathematical physics. The spectral theory of such problems has been investigated already by Birkhoff \cite{B08i, B08ii} and Tamarkin \cite{T28}.

The main result of this paper is Theorem \ref{maint}, where we establish the similarity of a general Dirac operator $L$ on an interval $[0,\omega]$ and an operator given by a direct sum of finite rank operators, all but one of which have rank at most two. A more explicit version of the result is in Theorem \ref{mainres11}. The similarity of the operators allows us to obtain structural results about $L$ and the group generated by the operator $iL$. 
In particular, the asymptotic estimates of the eigenvalues of $L$ appear in Theorem \ref{mainres1}. The (generalized) spectrality of operator $L$ is stated in Theorem \ref{mainresspec}. In Theorem
\ref{mainresequi}, we obtain a result on equiconvergence of the spectral decompositions in the Hilbert-Schmidt operator topology.  Finally, the group generated by the operator $iL$ is exhibited in Theorem \ref{mainres3}.

We study 
one-dimensional Dirac operators in one of their most general forms. We let $\mathcal{H}= L^2=L^2([0, \omega], \mathbb{C}^2)=L^2[0, \omega]\oplus L^2[0, \omega]$~--
be the Hilbert space (of equivalence classes) of $\mathbb{C}^2$-valued square-summable functions on $[0, \omega]$. The inner product in $\mathcal{H}$ is given by
$$
\la f, g\ra =\frac{1}{\omega}\int_0^\omega(f_1(t)\overline{g}_1(t)+f_2(t)\overline{g}_2(t))\,{d}t, \ f = (f_1, f_2), \ g=(g_1, g_2)\in \HH.
$$
The space $\HH$ is isometrically isomorphic to the space $L_{2, \omega}=L_{2, \omega}(\mathbb{R}, \CC^2)$ of $\omega$-periodic $\CC^2$-valued functions on $\mathbb{R}$ that are square-summable over $[0,\omega]$. In this paper, we typically do not distinguish these two spaces.

By $W_2^1([0, \omega], \mathbb{C}^2)$ we denote the Sobolev space of absolutely continuous $L^2$-functions with derivatives
in $L^2$ and the inner product
 $\langle f, g\rangle_W =\la f, g\ra+\la f', g'\ra$, $f, g\in W_2^1([0, \omega], \mathbb{C}^2)$.

We consider Dirac operators $L_{bc}: D(L_{bc})\subset \mathcal{H}\to\mathcal{H}$,
such that
\begin{equation}\label{dir1}
(L_{bc}y)(t)=i
\begin{pmatrix}
1 & 0 \\
0 & -1
\end{pmatrix}
\frac{dy}{dt} - P(t)y(t),
\end{equation}
where $t\in [0, \omega]$ and 
\begin{equation}\label{dir2}
P(t)=
\begin{pmatrix}
p_1(t) & p_2(t) \\
p_3(t) & p_4(t)
\end{pmatrix},
\end{equation}
$p_j\in L^2[0, \omega]$, $1\leqslant j\leqslant 4$.

The notation {\it ``bc''} refers to the various kinds of boundary conditions that we will use to describe the domain $D(L_{bc})$. We employ the notation from \cite{BDS11, DM06}
and consider the following three cases:
\ben
\item[(a)] periodic boundary condition ($bc = per$: $y(0)=y(\omega)\in\mathbb{C}^2$);
\item[(b)] anti-periodic boundary condition  ($bc=ap$: $y(0)=-y(\omega)\in\mathbb{C}^2$);
\item[(c)] Dirichlet boundary condition ($bc=dir$: $y_1(0)=y_2(0), y_1(\omega)=y_2(\omega)$), where $y=(y_1, y_2)\in W_2^1([0, \omega], \mathbb{C}^2)$.
\een
We let $D(L_{bc})=\{y\in W_2^1([0, \omega], \mathbb{C}^2), y\in bc\}$ and use the notation $L_{per}$, $L_{ap}$, or $L_{dir}$ to refer to the Dirac operators with periodic, anti-periodic or Dirichlet boundary conditions, respectively.

We remark that $per$ and $ap$ are Birkhoff regular but not strictly regular \cite{SS14}. The boundary conditions $dir$ are strictly regular. 
We also note an alternative form of Dirac operators commonly used in the literature \cite{SS15, SS14}:
\begin{equation}\label{dir2'}
(\bar{L}_{bc}u)(t)=
\begin{pmatrix}
0 & 1 \\
-1 & 0
\end{pmatrix}
\frac{du}{dt} -
\begin{pmatrix}
v_1(t) & v_2(t) \\
v_3(t) & v_4(t)
\end{pmatrix}
u(t),
\end{equation}
where $v_i\in L^2[0, \omega]$, $1\leqslant i\leqslant 4$. The study of $\bar{L}_{bc}$ is reduced to that of ${L}_{bc}$ via a simple change of variables: $u_1 = \frac12 (y_1 + y_2)$, $u_2 = \frac i2 (y_1 -y_2)$, $y = (y_1, y_2)$, $y \in L^2([0, \omega], \CC^2)$.

The research on Dirac operators has a long and illustrious history (see \cite{DM06, SS14, T92} and references therein). It was initiated in \cite{D28i, D28ii}, where Dirac considered operators of the form \eqref{dir2'} with $v_2=v_3=0$, $v_1=V-m$, $v_4=V+m$, where the function $V$ represented the potential of an electromagnetic field and $m$ represented mass. Since then, Dirac operators with various kinds of matrix potentials have become a staple in theoretical and applied research in mathematical physics. Multisection semiconductor lasers \cite{ZRS03} is one of the applications.
 In the next several paragraphs, we mention  theoretical papers that are most relevant to our research; the list is by no means exhaustive. There are also many papers that study operators that can be transformed into Dirac operators; \cite{BKR17, BKU18, RSS99} are among them.

In \cite{DM10, DM12JFA, DM12, DM06}, Djakov and Mityagin extensively studied the spectral theory of the operators $L_{bc}$ in the cases of periodic, anti-periodic and other Birkhoff regular boundary conditions. 
For example, in \cite{DM10}, generalized spectrality of Dirac operators was proved.

In \cite{LM14, LM15, MO12}, Malamud {\it et.~al.}~studied the questions of completeness of root vectors of Dirac operators.

In \cite{S16t, S16, SS15, SS14}, equiconvergence of spectral decompositions in the strong operator topology was obtained for different classes of Dirac operators in Hilbert and Banach spaces.

In \cite{BDS11}, the method of similar operators was used for the first time to study Dirac operators. Asymptotic formulas for the eigenvalues and equiconvergence of spectral decompositions in the uniform operator topology were obtained there for the case $p_1 = p_4 =0$ in \eqref{dir2}. The latter case is commonly considered by various authors \cite{DM06, S16, SS14}. 

In this paper, we continue the line of research started in \cite{BDS11} and consider the case of the general matrix potential. On the technical side, this necessitates a more elaborate version of
the method of similar operators, akin to the one developed in \cite{BP17}. On the results side, it turns out that a generic potential leads to a better estimate of the structure of the spectrum of the Dirac operator compared to the special case of \cite{BDS11}. 

 The remainder of the paper is organized as follows. In Section \ref{prelim}, we set the basic notation and exhibit a few preliminary results that will be used throughout the paper. We also state a version of our main result  -- Theorem \ref{maint} -- to serve as a goal for the following four sections. In Section \ref{sec3}, we describe the idea  of the method of similar operators and illustrate it with a two dimensional example. In Section \ref{sect3}, we present preliminary similarity transforms which put Dirac operators into a form that is amenable for the use of the method described in Section \ref{sec3}. The method calls for constructing two more similarity transforms, the first of which appears in Section \ref{transext} (see Theorem \ref{1trans}).
We remark that the approach we use to construct the first similarity transform is new and requires many technical details. In Section \ref{fsim}, we perform the last similarity transform of the method and finally arrive at the main result of the paper -- Theorem \ref{mainres11}.  We use Theorem \ref{mainres11} to obtain asymptotic estimates of the spectrum $\s(L_{bc})$ in Section \ref{sec7} and results on equiconvergence of spectral decompositions in Section \ref{sec8}. The final Section \ref{groupsec} is devoted to the description of the group generated by the operator $iL_{bc}$. The description is based on Theorem \ref{mainres11} and the asymptotic estimates of the spectrum from Section \ref{sec7}.

\section{Notation and Preliminaries}\label{prelim}

In this section, we introduce the necessary notation and state a version of the main result of the paper. A lot of the notation is the same as in \cite{BKU18}; the results, however, though similar in nature, apply to a very different class of operators.

We begin with notation that will allow us to formulate a technical condition on the potential matrix $P$ in \eqref{dir1} that we need to distinguish between different cases in our main results.

Given a function $v\in L^2[0, \omega]$, its Fourier series is 
$$
v(t)\sim\sum_{n\in\mathbb{Z}}\widehat{v}(n)e^{i\frac{2\pi n}{\omega}t}, \quad t\in [0, \omega],
$$
where the Fourier coefficients are
$$
\widehat{v}(n)=\frac{1}{\omega}\int_0^\omega v(t)e^{-i\frac{2\pi n}{\omega}t}\,dt=\la v, e_n\ra, \quad n\in\mathbb{Z}.
$$
As usual, we shall refer to the number $\widehat{v}(0)$ as the average (value) of the function $v$. Most commonly, we shall use the averages $\widehat{p}_1(0)$ and $\widehat{p}_4(0)$.
Moreover, we shall occasionally assume the following technical condition:
\begeq\label{badcond}
r = \frac{\omega\beta}{2\pi}=\frac\omega{2\pi}\left(\widehat{p}_1(0) - \widehat{p}_4(0)\right) \notin \ZZ\setminus\{0\}.
\eq
For certain results, we will have to use a different approach if the above condition does not hold.
We shall discuss it in more detail in Remark 
\ref{rembad2}.

To simplify the exposition, we shall also use the following notation: 
\begin{equation}\label{dir6}
\nu=\frac{1}{2}(\widehat{p}_1(0)+\widehat{p}_4(0)), \quad
\theta=-\pi r=\frac{\omega}{2}(\widehat{p}_4(0)-\widehat{p}_1(0));
\end{equation}
\begin{equation}\label{dir5}
\varphi(t)=\nu t-\int_0^tp_1(\tau)\,d\tau, \quad \psi(t)=-\nu t+\int_0^tp_4(\tau)\,d\tau,
\end{equation}
\begin{equation}\label{dir7}
q_2(t)=p_2(t)e^{i(\psi(t)-\varphi(t))}, \quad q_3(t)=p_3(t)e^{i(\varphi(t)-\psi(t))}, \quad t\in [0, \omega].
\end{equation}
We note that $\varphi(0)=\psi(0)=0$ and $\varphi(\omega)=\psi(\omega)=\theta$.

 We continue our exposition of definitions and notation with orthogonal direct sums of Hilbert spaces and operators.

Throughout this section and the following one, $\HH$ will denote an abstract Hilbert space. By $B(\HH)$ we shall mean the Banach algebra of all bounded linear operators in $\HH$.
 We shall also make use of  the ideal
of Hilbert-Schmidt operators in $\HH$ denoted by $\mathfrak{S}_2(\HH)$ and the ideal of nuclear operators -- $\mathfrak{S}_1(\HH)$. Recall that in $\mathfrak{S}_1(\HH)$ the norm is given by $\|X\|_1=\sum\limits_{n=1}^\infty |s_n|$, where $(s_n)$ is the
sequence of singular values of the operator $X$. The norm in $\mathfrak{S}_2(\HH)$  is
$\|X\|_2=(\mathrm{tr}\,XX^*)^\frac{1}{2} =\left(\sum\limits_{n=1}^\infty |s_n|^2\right)^{1/2}$. We refer to \cite{DS88II, GK65} for the standard properties of these ideals used in this paper.

Assume that  $\HH$ is an orthogonal direct sum of nontrivial closed subspaces $\mathcal{H}_n$, $n\in\mathbb{Z}$, i.e.
\begin{equation}\label{bask8}
\mathcal{H}=\bigoplus_{n\in\mathbb{Z}}\mathcal{H}_n,
\end{equation}
where $\mathcal{H}_m$ is orthogonal to $\mathcal{H}_n$ for $m\ne n$, $m$, $n\in\mathbb{Z}$, and $x=\sum\limits_{n\in\mathbb{Z}}x_n$,
$x_n\in\mathcal{H}_n$, $\|x\|^2=\sum\limits_{n\in\mathbb{Z}}\|x_n\|^2$. Such a representation of $\mathcal{H}$ generates a resolution of the identity
 $\{{P}_n$, $n\in\mathbb{Z}\}$, where the idempotents
${P}_n$, $n\in\mathbb{Z}$, have the following properties:
\ben
\item ${P}_n^*={P}_n$, $n\in\mathbb{Z}$.

\item ${P}_m{P}_n=0$ for $m\ne n$, $m$, $n\in\mathbb{Z}$.

\item The series $\sum\limits_{n\in\mathbb{Z}}{P}_nx$ converges unconditionally to $x\in\mathcal{H}$ and
$\|x\|^2=\sum\limits_{n\in\mathbb{Z}}\|{P}_nx\|^2$.


\item $\mathcal{H}_k=\mathrm{Im}\,{P}_k$, $x_k={P}_kx$, $k\in\mathbb{Z}$.
 \een

We remark that the third of the above properties is equivalent to
\ben
\item[3'.] ${P}_kx=0$ for all $k\in\mathbb{Z}$ implies $x=0$.
\een

 We note  that given a resolution of the identity $\{P_n: n\in\ZZ\}$ and $X\in\mathfrak S_2(\HH)$ we have
$$
\|X\|_2^2=\sum\limits_{m, n\in\ZZ}\|P_nXP_m\|_2^2.
$$

\begin{defn}\label{baskdef6}
We say that a linear operator ${A}: D({A})\subset\mathcal{H}\to\mathcal{H}$ is an \textit{orthogonal direct sum 
\begin{equation}\label{bask9}
{A}=\bigoplus_{n\in\mathbb{Z}}{A}_n,
\end{equation}
of bounded linear operators} ${A}_n\in B(\mathcal{H}_n)$, $n\in\mathbb{Z}$, with respect to a decomposition
(\ref{bask8}), 
if the following three conditions hold:

\ben
\item $D({A})=\{x\in\mathcal{H}: \sum\limits_{k\in\mathbb{Z}}\|{A}_kP_kx\|^2<\infty 
\}$ and $\mathcal{H}_n\subset D({A})$ for all  $n\in\mathbb{Z}$;

\item Each $\mathcal{H}_n$, $n\in\mathbb{Z}$, is an invariant subspace of the operator ${A}$ and
${A}_n$, $n\in\mathbb{Z}$, is the restriction of ${A}$ to $\mathcal{H}_n$, $n\in\mathbb{Z}$;

\item ${A}x=\sum\limits_{k\in\mathbb{Z}}{A}_kP_k x$, $x\in D({A})$.
\een
\end{defn}

We remark that for an orthogonal direct sum of operators we have $\sigma({A}_k)\subset\sigma({A})$, $k\in\mathbb{Z}$.
It may, however, happen that $\sigma({A})$ is strictly larger than the closure of the union of $\sigma({A}_k)$, $k\in\mathbb{Z}$. We also mention a slight abuse of notation in the above definition, where we treat the operators $A_k$ as both members of $B(\HH_k)$ and $B(\HH)$. We believe, it should be clear from the context what exactly is meant in each instance.

We shall illustrate the notion of direct sums of operators with the help of an operator
$L_{bc}^0$, which is the Dirac operator $L_{bc}$ with a trivial potential $P=0$, i.e. $p_i=0$, $1\leq  i\leq  4$. We shall call such operators \emph{unperturbed} or \emph{free}.
The operators $L_{bc}^0$, $bc\in\{per, ap, dir\}$, are self-adjoint and their spectral properties can be easily described \cite{BDS11, DM06}:
\ben
\item $\sigma(L_{per}^0)=\left\{\frac{2\pi n}{\omega}\right\}$, $n\in\mathbb{Z}$, where each $\lambda_n=\frac{2\pi n}{\omega}$, 
$n\in\mathbb{Z}$, is an eigenvalue of multiplicity two. Moreover, the corresponding eigenspace is given by $\HH_n = E_n^0=\mathrm{span}\,\{e_n^1, e_n^2\}$, where
$$
e_n^1=
\begin{pmatrix}
e_{-n} \\
0
\end{pmatrix},
\quad
e_n^2 =
\begin{pmatrix}
0 \\
e_n
\end{pmatrix},
\quad e_n(t)
=e^{i\lambda_nt}, \quad t\in [0, \omega].
$$

\item $\sigma(L_{ap}^0)=\left\{\frac{\pi(2n+1)}{\omega}\right\}$, $n\in\mathbb{Z}$, where each $\lambda_n=\frac{\pi(2n+1)}{\omega}$,
$n\in\mathbb{Z}$, is again an eigenvalue of multiplicity two and the corresponding eigenspace
is defined the same way as for $L^0_{per}$.

\item $\sigma(L_{dir}^0)=\left\{\frac{\pi n}{\omega}\right\}$, $n\in\mathbb{Z}$, where each
$\lambda_n=\frac{\pi n}{\omega}$, $n\in\mathbb{Z}$, is a simple eigenvalue, and the corresponding normalized eigenfunction is given by
 $s_n=\frac{1}{\sqrt{2}}(e_n^1+e_n^2)$, $n\in\mathbb{Z}$.
\een

We remark that Birkhoff regularity of the boundary conditions yields spectrality of the unperturbed operator $L^0_{bc}$ in the sense of Dunford \cite{DS88III}.

\bex\label{direx}
The Riesz projections of an operator $L_{bc}^0$ form resolutions of the identity that lead to the direct sums of operators. In particular, we denote by $P_n$ the (spectral) Riesz projection that corresponds to the eigenvalue $\l_n$ of the operator $L_{bc}^0$, $n\in\ZZ$. 
Then $\{P_n, n\in\ZZ\}$ is indeed a resolution of the identity that satisfies the three properties preceding Definition \ref{baskdef6}.
We shall also make use of the spectral projections $P_{(k)} = \sum\limits_{|n|\le k} P_n$,  and the corresponding coarser resolutions of the identity given by $\{P_{(k)}\}\cup\{P_n, |n|> k\}$, $k\in \ZZ_+ = \NN\cup\{0\}$. We then let $\HH_n = {\rm{Im}}\,P_n$ and $\HH_{(k)} = {\rm{Im}}\,P_{(k)}$, $n, k\in\ZZ$, $k\ge 0$, be the corresponding eigenspaces, which provide us with orthogonal decompositions of $\HH = L^2$. The operators
$(L_{bc}^0)_n = L_{bc}^0\vert_{\HH_n} =\l_n I_n$, where $I_n$ is the identity operator in $\HH_n$, and $(L_{bc}^0)_{(k)}= L_{bc}^0\vert_{\HH_{(k)}}$ 
then give orthogonal direct sums of the operator $L_{bc}^0$:
$$
L_{bc}^0=(L_{bc}^0)_{(k)}\oplus\left(\bigoplus_{|n|>k}(L_{bc}^0)_n\right)=(L_{bc}^0)_{(k)}\oplus
\left(\bigoplus_{|n|>k}\lambda_nI_n\right), \quad k\in\mathbb{Z}_+,
$$
with respect to the orthogonal decompositions of $\mathcal{H}$ given by
$$
\mathcal{H}=\mathcal{H}_{(k)}\oplus\left(\bigoplus_{|n|>k}\mathcal{H}_n\right), \quad k\in\mathbb{Z}_+.
$$
We remark that the operators $(L_{bc}^0)_n$ have rank $2$ if $bc\in\{per, ap\}$ and rank $1$ if $bc = dir$.
\eex

To state our main result we need the following standard notion of similarity between unbounded linear operators.

\bd\label{defsim}
Two linear operators $A_m: D(A_m)\subset\mathcal{H}\to\mathcal{H}$, $m=1, 2$, are called 
\emph{similar}, if there exists a continuously invertible operator
 $U\in B(\mathcal{H})$ such that
$$
A_1Ux=UA_2x, \quad x\in D(A_2), \quad UD(A_2)=D(A_1).
$$
The operator $U$ is called the \emph{similarity transform} of $A_1$ into $A_2$.
\ed

The following two definitions are natural in view of the notion of similarity.

\bd\label{baskdef7'}
Assume that $\HH$ has an orthogonal decomposition \eqref{bask8}.
Given an invertible operator $U\in B(\HH)$,
we call the direct sum
\begin{equation}\label{bask10''}
\mathcal{H}=\bigoplus_{k\in\mathbb{Z}}U\mathcal{H}_k,
\end{equation}
a quasi- or $U$-orthogonal decomposition. If $U=I+W$ for some operator $W\in\mathfrak{S}_2(\mathcal{H})$, we call such a decomposition of  $\mathcal{H}$ a \emph{Riesz decomposition}.
\ed
 
 We remark that a $U$-orthogonal decomposition (\ref{bask10''}) of $\mathcal{H}$ may be regarded as an orthogonal decomposition of $\mathcal{H}$ with respect to an equivalent inner product
$$
\la x, y\ra_U=\la Ux, Uy\ra, \quad x, y\in\mathcal{H}.
$$

\bd\label{baskdef8'}
Assume that a linear operator ${A}: D({A})\subset\mathcal{H}\to\mathcal{H}$ is an orthogonal direct sum of the form \eqref{bask9}. Assume also that $\widetilde A$ is similar to $A$ and $U\in B(\HH)$ is the similarity transform of $\widetilde A$ into $ A$.
We then say that $\widetilde A$ is a \emph{quasi- or $U$-orthogonal direct sum} 
$$
 \widetilde A=\bigoplus_{k\in\mathbb{Z}}\widetilde{ A}_k
$$
of bounded linear operators $\widetilde{ A}_k$, $k\in\mathbb{Z}$, with respect to a decomposition of the space
$\mathcal{H}$ of the form (\ref{bask10''}), if $\widetilde{ A}_k=U A_kU^{-1}$, $k\in\mathbb{Z}$.
\ed

Directly from Definition \ref{defsim}, we have the following result about the spectral properties of similar operators.

\bl\label{basklh1}
Let $A_m: D(A_m) \subset \HH\to\HH$, $m = 1,2$, be two similar operators with the operator $U$ being the similarity transform of $A_1$ into $A_2$. Then the following properties hold.
\ben
\item[(1)]	We have $\s(A_1) = \s(A_2)$, $\s_p(A_1) = \s_p(A_2)$, and $\s_c(A_1) = \s_c(A_2)$, where $\s_p$ denotes the point spectrum and $\s_c$ denotes the continuous spectrum;
\item[(2)] If $\l$ is an eigenvalue of the operator $A_2$ and $x$ is a corresponding eigenvector, then $y = Ux$ is  an eigenvector of the operator $A_1$ corresponding to the same eigenvalue $\l$.
\item[(3)]	Assume that the operator $A_2$ is an  orthogonal direct sum $A_2 = \bigoplus\limits_{n\in\ZZ} (A_{2})_n$ with respect to an orthogonal decomposition $\HH = \bigoplus\limits_{n\in\ZZ} \HH_n$. 
Then the operator $A_1$ is a $U$-orthogonal direct sum $A_{1} = \bigoplus\limits_{n\in\ZZ} (A_{1})_n$ with respect to the $U$-orthogonal  decomposition $\HH = \bigoplus\limits_{n\in\ZZ} \widetilde{\HH}_n$, where ${\widetilde{\HH}_n} = U\HH_n$. Moreover, if $\{P_n\}$ is the resolution of the identity corresponding to the decomposition $\HH = \bigoplus\limits_{n\in\ZZ} \HH_n$, then $\{\widetilde P_n = UP_nU^{-1}\}$ is the  resolution of the identity corresponding to the decomposition $\HH = \bigoplus\limits_{n\in\ZZ} \widetilde{\HH}_n$.
\item[(4)] If $A_2$ is a generator of a $C_0$-semigroup (or group) $T_2: \mathbb{J}\to 
B(\mathcal{H})$, $\mathbb{J}\in\{\RR,\RR_+\}$, then the operator $A_1$ generates the $C_0$-semigroup (group)
$$
T_1(t)=UT_2(t)U^{-1}, \quad t\in\mathbb{J}, \quad T_1: \mathbb{J}\to B(\mathcal{H}),\  \mathbb{J}\in\{\RR,\RR_+\}.
$$
\een
\el

The following is the key result of the paper. In its formulation, we use the notation introduced when we represented the operator $L_{bc}^0$ as the direct sums of bounded operators.

\bt\label{maint}
Assume that an operator $L_{bc}$ is given by \eqref{dir1}. 
There exists $m\in\mathbb{Z}_+$, such that the operator  $L_{bc}$ is similar to the operator $\widetilde{L}_{bc}^P-V$, where 
 $V\in\mathfrak{S}_2(\mathcal{H})$, and the subspaces $\mathcal{H}_{(m)}$ and $\mathcal{H}_n$, $|n|>m$, are invariant for the operators $\widetilde{L}_{bc}^P$ and 
 $V$. Moreover, the dimension of $\HH_{(m)}$ is at most $2m+2$, the dimensions of $\mathcal{H}_n$, $|n|>m$, are at most $2$, the operator $\widetilde{L}_{bc}^P$ can be written explicitly, we have
$$
L_{bc}W_{bc}(I+U)=W_{bc}(I+U)(\widetilde{L}_{bc}^P-V),
$$
and the operator $L_{bc}$ is a $W_{bc}(I+U)$-direct sum
$$
L_{bc}=W_{bc}(I+U)\left(\left(\widetilde{L}_{bc}^P-V\right)_{(m)}\oplus\left(\bigoplus_{|n|>m}\left(\widetilde{L}_{bc}^P-V\right)_n\right)\right)(I+U)^{-1}W_{bc}^{-1},
$$
for some $W_{bc}\in B(\HH)$ and $U\in\mathfrak{S}_2(\mathcal{H})$.
\et

A more explicit version of the above result is in Theorem \ref{mainres11}.


To describe the $C_0$-group generated by the operator $iL_{bc}$ in Section \ref{groupsec}, we shall make use of the following lemma.

\bl\label{baskuskova_lh6}\cite[Lemma 3.4]{BKU18}
Assume that an operator $ A$ is an orthogonal direct sum as in Definition \ref{baskdef6}. Then $A$ is a generator of a $C_0$-group of operators
$T: \mathbb{R}\to B(\mathcal{H})$ if and only if
\begin{equation}\label{bask_16'}
\sup_{|t|\leq b}\sup_{n\in\mathbb{Z}}\|e^{t A_{n}}\|_{B(\mathcal{H}_n)}=C(b)<\infty,
\end{equation}
 $b\ge  1$. If (\ref{bask_16'}) holds, then the operators $ T(t)$, $t\in\mathbb{R}$, are
 orthogonal direct sums
$$
 T(t)=\bigoplus_{n\in\mathbb{Z}}e^{t A_{n}}, \quad t\in\mathbb{R},
$$
with respect to the orthogonal decomposition  (\ref{bask8}) of $\mathcal{H}$.
\el

%


\section{The method of similar operators}\label{sec3}
The method of similar operators has its origins in various similarity and perturbation techniques. Among them, we mention classical perturbation methods of celestial mechanics, Ljapunov's kinematic
similarity method \cite{GKK96, Lj56, N15},  Friedrichs' method of similar operators that is used in quantum mechanics \cite{F65}, and Turner's method of similar operators  \cite{T65, U04}. A close relative of the method is the Krylov-Bogolyubov substitution described in \cite{B84u}.

The method of similar operators has been extensively developed and used for various classes of unbounded linear operators, see e.g. \cite{B83, B85, B86, B94, B99, B15, BDS11, BK88, BKR17, BKU18}. In this paper, we use a version of the method that is derived mostly from \cite{BDS11, BP17}. In Subsection \ref{basmeth}, we exhibit the basic ideas and theorems of the method. In the following Subsection \ref{genJG}, we first illustrate the method in the case of $2\times 2$ matrices and next provide a construction for direct sums of operators.

\subsection{The idea and basic theorems of the method of similar operators.}\label{basmeth}

The main idea of the method is to construct a similarity transform for an operator $A - B:  D(A) \subset \HH \to \HH$, where the spectrum of the operator $A$ is known and has certain properties, and the operator $B$ is $A$-bounded (see Definition \ref{abdd} below). The goal of the method is to obtain an operator $B_1$ such that the operator $A-B$ is similar to $A-B_1$ and  the spectral properties of $A-B_1$ are in some sense close to those of $A$. In particular, certain spectral subspaces of $A$ are mapped by the similarity transform onto certain subspaces that are invariant for $A-B_1$.
In this paper, the role of $A$ is played by the operator $W_{bc}^{-1}L_{bc}^0W_{bc}$ (see the notation in Theorem \ref{maint} and Section \ref{sect3}  below).

\bd\label{abdd}
 Let $A: D(A) \subset \HH \to \HH$ be a linear operator. A linear operator $B: D(B) \subset \HH \to \HH$ is $A$-bounded if $D(B) \supseteq D(A)$ and $\|B\|_A = \inf\{c > 0: \|Bx\| \le c(\|x\| + \|Ax\|),\ x \in D(A)\} < \infty$.
\ed
The space $\mathfrak L_A(\HH)$ of all $A$-bounded linear operators is a Banach space with respect to the norm $\|\cdot\|_A$. Moreover, given $\l_0 \in \rho(A)$, where $\rho(A) = \CC\backslash\s(A)$ is the resolvent set of $A$, we have $B \in \mathfrak L_A(\HH)$ if and
only if $B(\l_0I - A)^{-1} \in B(\HH)$ and $\|B\|_{\l_0} = \|B(\l_0I - A)^{-1}\|_{B(\HH)}$ defines an equivalent norm in $\mathfrak L_A(\HH)$ \cite{EN00}.

The method of similar operators uses the \emph{commutator transform} 
$\mathrm{ad}_A:
D(\mathrm{ad}_A)\subset B(\mathcal{H})\to B(\mathcal{H})$ defined by
\begeq\label{ct}
\mathrm{ad}_AX = AX-XA, \quad X\in D(\mathrm{ad}_A).
\eq
The domain $D(\mathrm{ad}_A)$ in \eqref{ct} consists of all  $X\in B(\mathcal{H})$ such that
the following two properties hold:
\ben
\item $XD(A)\subseteq D(A)$;
\item The operator $\mathrm{ad}_AX: D(A)\to\mathcal{H}$ admits a unique extension to a bounded operator $Y\in B(\mathcal{H})$; we then let
$\mathrm{ad}_AX=Y$.
\een

The key notion of the method of similar operators is that of an admissible triplet. Once such a triplet is constructed, achieving the goal of the method becomes a  routine task.

\bd[\cite{BDS11, BP17}]\label{baskdef8}
Let $\mathcal{M}$ be a linear subspace of $\mathfrak{L}_A(\mathcal{H})$,
$J: \mathcal{M}\to\mathcal{M}$, and $\Gamma: \mathcal{M}\to B(\mathcal{H})$.
The collection $(\mathcal{M}, J, \Gamma)$ is  an \emph{admissible triplet} for the operator $A$, and the space
$\mathcal{M}$ is the \emph{space of admissible perturbations}, if the following six properties hold.

\ben
\item\label{adprop1}
 $\mathcal{M}$ is a Banach space that is continuously embedded in $\mathfrak{L}_A(\mathcal{H})$, i.e., $\mathcal{M}$ has a norm $\|\cdot\|_\ast$
such that there is a constant $C>0$ that yields $\|X\|_A\le  C\|X\|_\ast$ for any
$X\in\mathcal{M}$.

\item\label{adprop2}
 $J$ and $\Gamma$ are bounded linear operators; moreover, $J$ is an idempotent.

\item\label{keyprop} $(\Gamma X)D(A)\subset D(A)$ and
$$
(\mathrm{ad}_A\,\Gamma X)x = (X-JX)x, \quad x\in D(A), \quad   X\in\mathcal{M};
$$
moreover $Y = \Gamma X\in B(\mathcal{H})$ is the unique solution of the equation
\begin{equation}\label{bask11'}
\mathrm{ad}_A\,Y = AY-YA = X-JX,
\end{equation}
that satisfies $JY=0$.

\item\label{adprop4} $X\Gamma Y$, $(\Gamma X)Y\in\mathcal{M}$ for all $X, Y\in\mathcal{M}$, and there is a constant $\gamma>0$ such that
$$
\|\Gamma\|\le\gamma, \quad \max\{\|X\Gamma Y\|_\ast, \|(\Gamma X)Y\|_\ast\}\le \gamma\|X\|_\ast\|Y\|_\ast.
$$

\item\label{adprop5} $J((\Gamma X)JY)=0$ for all $X, Y\in\mathcal{M}$.

\item\label{adprop6} For every $X\in\mathcal{M}$ and $\varepsilon>0$ there exists a number  $\lambda_\varepsilon\in\rho(A)$,
such that $\|X(A-\lambda_\varepsilon I)^{-1}\|<\varepsilon$.
\een
\ed

To illustrate  the above definition, one should think of the operators involved in terms of infinite matrices. The operator $A$ is then represented by an infinite diagonal matrix and the operator $B$ -- by a matrix with some kind of off-diagonal decay. The transform $J$ should be thought of as a projection that picks  the main (block) diagonal of an infinite matrix, whereas the transform $\Gamma$ annihilates the main (block) diagonal and weighs the remaining diagonals in accordance with equation \eqref{bask11'} thereby introducing or enhancing the off-diagonal decay. A more precise illustration is provided in Section \ref{genJG}.

To formulate the main theorem of the method of similar operators for an operator $A-B$, we use the function
$\Phi:\M\to\M$ given by
\begeq\label{bask13}
\Phi(X) = B\Gamma X-(\Gamma X)(JB)-(\Gamma X)J(B\Gamma X)+B.
\eq 
\bt[\cite{BDS11, BP17}]\label{baskth6}
Assume that $(\mathcal{M}, J, \Gamma)$ is an admissible triplet for an operator $A: D(A)\subset\mathcal{H}\to\mathcal{H}$ 
and $B\in\mathcal{M}$. Assume also that
\begeq\label{bask12}
4\gamma\|J\|\|B\|_\ast<1,
\eq
where $\gamma$ comes from the Property 4 of Definition \ref{baskdef8}. Then 
the map $\Phi:\mathcal{M}\to\mathcal{M}$ given by \eqref{bask13} is a contraction and has a unique fixed point 
$X_*$ in the ball 
\begeq\label{ball}
\{X\in\mathcal{M}: \|X-B\|_\ast\le  3\|B\|_\ast\}, 
\eq
which can be found as a limit of simple iterations: $X_0=0$, $X_1=\Phi(X_0) = B$, etc.
Moreover,  the 
operator $A-B$ is similar to the operator $A-JX_*$ 
and the similarity transform of 
$A-B$ into $A-JX_*$ is given by $I+\Gamma X_*\in
B(\mathcal{H})$. 
\et

The space $\mathcal M$ in the above theorem is typically constructed based on the  properties of the operator $B$. Condition \eqref{bask12} is there to guarantee existence
of the solution of the functional equation $\Phi(X) = X$ or, in other words, existence and uniqueness of the fixed point of $\Phi$ in a ball in the space $\mathcal M$. 

We will need the following consequence of Lemma \ref{basklh1} and Theorem \ref{baskth6}.
\bt[\cite{BP17}]\label{baskth7}
Assume that $(\mathcal{M}, J, \Gamma)$ is an admissible triplet for $A: D(A)\subset\mathcal{H}\to\mathcal{H}$, 
$B\in\mathcal{M}$ satisfies (\ref{bask12}), and $A-JX_*$ is a generator of a $C_0$-group $\widetilde{T}: \mathbb{R}\to B(\mathcal{H})$. Then the operator $A-B$ is a generator of the $C_0$-group $T: \mathbb{R}\to\mathrm{End}\,\mathcal{H}$ given by
$$
T(t)=(I+\Gamma X_*)\widetilde{T}(t)(I+\Gamma X_*)^{-1}, \quad t\in\RR,
$$
where $X_*$ is the fixed point of the function $\Phi$ in (\ref{bask13}).
\et

In many cases, it can be difficult to define the space $\mathcal{M}$ of admissible perturbations for a given operator $A-B$. It may, however, be possible to pick a good space $\mathcal{M}$ first, and then find an operator $A-C$ that is similar to $A-B$ and such that $C\in \mathcal{M}$. In our case, this will be possible if the following assumption holds.

\begin{ass}[\cite{BP17}]\label{baskpred1}
Assume that $(\mathcal{M}, J, \Gamma)$ is  an {admissible triplet} for an operator $A$
such that the transforms $J$ and $\Gamma$ are restrictions of linear operators from 
 $\mathfrak{L}_A(\mathcal{H})$ to  $\mathfrak{L}_A(\mathcal{H})$ denoted by the same symbols. Assume also that the operator $B\in \mathfrak{L}_A(\mathcal{H})$ has the following five properties.
 
\ben
\item\label{propas1} $\Gamma B\in B(\mathcal{H})$ and $\|\Gamma B\|<1$;

\item\label{propas2} $(\Gamma B)D(A)\subset D(A)$;

\item\label{mprop} $B\Gamma B$, $(\Gamma B)JB\in\mathcal{M}$;

\item\label{propas4} $A(\Gamma B)x-(\Gamma B)Ax=Bx-(JB)x$, $x\in D(A)$;

\item\label{propas5} For any  $\varepsilon>0$ there is $\lambda_\varepsilon\in\rho(A)$ such that
$\|B(A-\lambda_{\varepsilon}I)^{-1}\|<\varepsilon$.
\een
\end{ass}
\bt[\cite{BP17}]\label{baskth8}
If Assumption \ref{baskpred1} holds then the operator $A-B$ is similar to   $A-JB-B_0$, where
$B_0=(I+\Gamma B)^{-1}(B\Gamma B-(\Gamma B)JB)$. 
The similarity transform is given by $I+\Gamma B$ so that
$$
(A-B)(I+\Gamma B)=(I+\Gamma B)(A-JB-B_0).
$$
\et


\subsection{Constructing the transforms $J$ and $\Gamma$ for direct sums of operators.}\label{genJG}

We begin with the following, simplest possible non-trivial illustration of the method of similar operators that will prove surprisingly useful for us.

\bex\label{twotwo}
Let $\HH = \CC^2$, so that $B(\HH)$ is simply the set of all $2\times2$ matrices with the operator norm, and $\mathfrak S_2(\HH)$ is the same set endowed with the Frobenius norm.
We let 
$$A = \left(
\begin{array}{cc}
a_{11} & 0 \\
0 & a_{22} 
\end{array}
\right)
\quad\mbox{and}\quad
B = \left(
\begin{array}{cc}
b_{11} & b_{12} \\
b_{21} & b_{22} 
\end{array}
\right),
$$
where $a_{11}, a_{22}, 
b_{11}, b_{12}, 
b_{21}, b_{22} \in \CC$ and $a_{11}\neq a_{22}$.
With $\mathcal M \in \{B(\HH), \mathfrak S_2(\HH)\}$, we define
\[
J\left(
\begin{array}{cc}
x_{11} & x_{12} \\
x_{21} & x_{22} 
\end{array}
\right) = \left(
\begin{array}{cc}
x_{11} & 0 \\
0 & x_{22} 
\end{array}
\right)
\mbox{ and }
\]
\[\Gamma \left(
\begin{array}{cc}
x_{11} & x_{12} \\
x_{21} & x_{22} 
\end{array}
\right) =\left(
\begin{array}{cc}
0 & \frac{x_{12}}{a_{11}-a_{22}} \\
\frac{x_{21}}{a_{22}-a_{11}} & 0 
\end{array}
\right),\ 
\left(
\begin{array}{cc}
x_{11} & x_{12} \\
x_{21} & x_{22} 
\end{array}
\right)
\in\mathcal M.
\]
\eex
One can easily verify that $\{\mathcal M, J, \Gamma\}$ is then an admissible triplet with $\|J\| = 1$ and $\gamma = |a_{22}-a_{11}|^{-1}$ (see Definition \ref{baskdef8}). Therefore, Theorem 
\ref{baskth6} is applicable for $A-B$ as long as $\|B\|_* < \frac14|a_{11}-a_{22}|$. In that case, the operator $A-B$ is similar to $A-JX_*$, where $X_*$ is a fixed point of \eqref{bask13}. In view of Property \ref{adprop5} of Definition \ref{baskdef8}, we have $JX_* = J(B\Gamma X_*) +JB$. Denoting 
\begeq\label{bgx}
X_* = \left(
\begin{array}{cc}
x^*_{11} & x^*_{12} \\
x^*_{21} & x^*_{22} 
\end{array}
\right)
\mbox{ and }
B\Gamma X_* =
\left(
\begin{array}{cc}
c_{11} & c_{12} \\
c_{21} & c_{22} 
\end{array}
\right),
\eq 
we then have that $A-B$ is similar to
\[
\left(
\begin{array}{cc}
a_{11}-b_{11}-c_{11} & 0 \\
0 & a_{22}-b_{22}-c_{22} 
\end{array}
\right)
\]
with the similarity transform given by
\[
I+\Gamma X_* = \left(
\begin{array}{cc}
1 & \frac{x^*_{12}}{a_{11}-a_{22}} \\
\frac{x^*_{21}}{a_{22}-a_{11}} & 1 
\end{array}
\right).
\]

\brem
In applying the method of similar operators to a Dirac operator, we will often end up with a sequence of $2\times2$ matrices. The usefulness of the above example will then be revealed.
\erem

Next, we generalize the construction of $J$ and $\Gamma$ from Example \ref{twotwo} to the case when an operator $A$ is an orthogonal direct sum induced by a spectral decomposition. We assume that the  resolvent  operator $(A-\l I)^{-1}$ is in $\mathfrak S_2(\HH)$ for any $\l \in \rho(A)$. 
 Additionally, we assume that the spectrum $\s(A)$ is a separated set: 
\begeq\label{speccond}
\s(A) = \bigcup_{n\in\ZZ} \{\l_n\} \mbox{ with } \delta =\inf\{|\l_m -\l_n|: m\neq n\in\ZZ\} > 0.
\eq
By $P_n$ we denote the (orthogonal) Riesz projection that corresponds to the spectral component $\{\l_n\}$, $n\in\ZZ$. We then have that the family $\{P_n, n\in\ZZ\}$ forms a resolution of the identity and the subspaces $\HH_n = {\rm Im } P_n$ provide an orthogonal decomposition of $\HH$ of the form \eqref{bask8}. As in Example \ref{direx}, we get $A = \bigoplus\limits_{n\in\ZZ} \l_n I_n$, where $I_n$ is the identity operator on $\HH_n$, $n\in\ZZ$.

Thus, an operator $X\in D(X)\subseteq \HH\to \HH$ is determined by an operator matrix with entries $X_{mn}=P_mXP_n$, $m,n\in\ZZ$, as long as the domain $D(X)$ is properly specified. As we mentioned above, the transform $J$  is supposed to pick out the main diagonal of the operator matrix. Hence, we define $J$ via
\begeq\label{tJ2}
(JX)_{mn} = \delta_{m-n} X_{mn},\ m,n\in\ZZ,\  X: D(X)\subseteq \HH\to \HH,
\eq
where $\delta_{k}$ is the usual Kronecker delta. Clearly, we have
\begeq\label{tJ}
JX = \sum_{n\in \ZZ} P_n XP_n, \ X\in \mathfrak S_2(\HH), 
\eq
where the series converges unconditionally in $\mathfrak S_2(\HH)$. 

Observe that in this setting the matrix of the commutator  $\mathrm{ad}_AX$ in \eqref{ct} satisfies
\[
(\mathrm{ad}_AX)_{mn} =   (\l_m - \l_n)X_{mn}, \ X\in D(\mathrm{ad}_A).
\]
Therefore, it is natural to define the transform $\Gamma$  via
\begeq\label{tG2}
(\Gamma X)_{mn} =
\begin{cases}
 \frac{1}{\l_m-\l_n}X_{mn}, & m\neq n;\\
0, & m=n; 
\end{cases}
\quad X: D(X)\subseteq\HH\to \HH. 
\eq
For $X\in \mathfrak S_2(\HH)$, we have
\begeq\label{tG}
\Gamma X =\sum_{\substack{m,n \in\ZZ \\ {m\neq n}}}\frac{P_mXP_n}{\l_m-\l_n}, %
\eq
where the series 
converges unconditionally in $\mathfrak S_2(\HH)$ because of \eqref{speccond}. 
In \eqref{tJ2} and \eqref{tG2}, we let the domains of the operators $JX$ and $\Gamma X$   be the largest possible, unless otherwise specified.
We will provide a more careful treatment of the case of $X \in \mathcal L_A(\HH)$ in  Section \ref{transext}.


\bl\label{JG1}
Assume that $\M = \mathfrak S_2(\HH)$, the transforms $J$ and $\Gamma$ are defined by \eqref{tJ} and \eqref{tG}, and for any $\varepsilon > 0$ there is $\lambda_\varepsilon\in\rho(A)$,
such that $\|(A-\lambda_\varepsilon I)^{-1}\|<\varepsilon$.
Then $(\M, J, \Gamma)$ is an admissible triplet.
\el

\bpf
Observe that Properties \ref{adprop1}, \ref{adprop2}, and \ref{adprop6} of Definition \ref{baskdef8} are automatically satisfied. Property \ref{adprop5} follows easily from \eqref{tJ} and \eqref{tG} via a straightforward computation. It can also be obtained as a direct application of \cite[Corollary 7.8]{BK05}. Property \ref{adprop4} clearly follows from \eqref{speccond} by letting $\gamma = \delta^{-1}$.

To prove Property \ref{keyprop} of Definition \ref{baskdef8}, 
 we first observe that $\Gamma \mathfrak S_2(\HH) \subseteq \mathfrak S_2(\HH)$ because of \eqref{speccond}.

Secondly, we need to prove that $(\Gamma X)D(A)\subset D(A)$ for $X \in \mathfrak S_2(\HH)$. Pick $x \in D(A)$ and $\l \in \rho(A)$. Then there is $y\in\HH$ such that 
$$
x = (A-\l I)^{-1}y = \sum_{n\in\ZZ} \frac1{\l_n-\l}P_ny.
$$
We then have
\[
\bs
\Gamma Xx & = (\Gamma X)(A-\l I)^{-1}y = \sum_{\substack{m,n \in\ZZ \\ {m\neq n}}}\frac{P_mXP_n y}{(\l_m-\l_n)(\l_n-\l)} \\ & 
= \sum_{\substack{m,n \in\ZZ \\ {m\neq n}}}\frac{P_mXP_n y}{(\l_m-\l_n)(\l_m-\l)} +
\sum_{\substack{m,n \in\ZZ \\ {m\neq n}}}\frac{P_mXP_n y}{(\l_m-\l)(\l_n-\l)} \\ &
= (A-\l I)^{-1}(\Gamma X)y +(A-\l I)^{-1}(\Gamma X)(A-\l I)^{-1}y \\ & =
 (A-\l I)^{-1}(\Gamma X)(x+y) \in D(A).
\end{split}
\]
Thirdly, for $X\in \mathfrak S_2(\HH)$ and $x\in D(A)$, we have
\[
\bs
A(\Gamma X)x - (\Gamma X)Ax &= \sum_{\substack{m,n \in\ZZ \\ {m\neq n}}}\frac{\l_mP_mXP_nx}{\l_m-\l_n} -
\sum_{\substack{m,n \in\ZZ \\ {m\neq n}}}\frac{\l_n P_mXP_nx}{\l_m-\l_n} \\ &= \sum_{\substack{m,n \in\ZZ \\ {m\neq n}}}{P_mXP_n}x =  (X - JX)x.
\end{split}
\]
Finally, we observe that $J(\Gamma X) = 0$ for all $X\in \mathfrak S_2(\HH)$.
\epf

As in Example \ref{direx}, we will also need to consider coarser resolutions of the identity given by  $\{P_{(k)}\}\cup\{P_n, |n|> k\}$, $k\in \ZZ_+$, where $P_{(k)} = \sum\limits_{|n|\le k} P_n$.
We then have two families of transforms  -- $\{J_k$, $k\in\ZZ_+\}$ and $\{\Gamma_k$, $k\in\ZZ_+\}$ --  that are given by 
\begeq\label{Jk}
J_kX = JX - P_{(k)}(JX)P_{(k)}+ P_{(k)}XP_{(k)}= P_{(k)}XP_{(k)}+\sum_{|n|>k} P_n XP_n, 
\eq
and
\begeq\label{Gk}
\Gamma_k X = \Gamma X- P_{(k)}(\Gamma X)P_{(k)}=\sum_{\substack{\max\{|m|,|n|\}> k \\ {m\neq n}}}\frac{P_mXP_n}{\l_m-\l_n}, \quad X\in\mathfrak S_2(\HH). %
\eq
 Clearly, $J_0$ and $\Gamma_0$ coincide with the transforms given by
\eqref{tJ} and \eqref{tG}, respectively, and 
$J_k X$ and $\Gamma_k X$, $k\in\ZZ_+$, are finite-rank perturbations of $JX = J_0X$ and $\Gamma X = \Gamma_0 X$, $X \in \mathfrak S_2(\HH)$. Moreover, 
\begeq\label{limG}
\lim_{k\to \infty}\Gamma_k X = 0
\eq
in the topology of $\mathfrak S_2(\HH)$.
Similarly to Lemma \ref{JG1}, one can check that $(\mathfrak S_2(\HH), J_k, \Gamma_k)$ is an admissible triplet for each $k\in\ZZ_+$ as long as the condition on the resolvent of $A$ holds. 

In Section \ref{transext}, we will show that in the case of Dirac operators the above construction will lead to an admissible triplet that satisfies Assumption \ref{baskpred1}.

\section{Preliminary similarity transforms for  Dirac operators}\label{sect3}

As indicated by the discussion at the end of Section \ref{basmeth}, one often cannot directly apply Theorem  \ref{baskth7} without using the transform from Theorem \ref{baskth8}.
In the case of a Dirac operator with a general potential, the situation is even more complicated. Before we can extend the construction of the transforms $J$ and $\Gamma$ from Section \ref{genJG} for Theorem \ref{baskth8}, we need yet another  preliminary similarity transform. This preliminary transform will depend on the type of boundary conditions imposed on the operator.

\brem\label{redrem}
In the literature (see e.g. \cite{DM06, S16, SS15, SS14}), the operators $L_{bc}$ are typically considered in the case when $p_1 = p_4 =0$ in \eqref{dir2}. For this case, as well as for the diagonal case $p_2 = p_3 =0$, an extensive spectral theory of Dirac operators have been created. It is often said \cite{S16, SS15, SS14} that the general case can be reduced to that of 
$p_1 = p_4 =0$ via a kind of similarity transform that is not connected to the method of similar
operators. Such a reduction, however, may change the spectral properties of the unperturbed
operator. This leads to complications that prevent one from clear understanding of the spectral theory of the general case.
\erem

\subsection{The preliminary similarity transform for $L_{dir}$.}

In this case, we can use the similarity transform from \cite{SS15, SS14} that was alluded to in Remark \ref{redrem}. We include the transform in the statement of the following theorem, which is merely a reformulation of the results in  \cite{SS15, SS14}.

\bt\label{dirth4}
An operator  $L_{dir}$ is similar to the operator $\widetilde{L}_{dir}: D(L_{dir})\subset\HH\to \HH$, given by
$$
(\widetilde{L}_{dir}y)(t)=i
\begin{pmatrix}
1 & 0 \\
0 & -1
\end{pmatrix} \frac{dy}{dt}
-\nu y(t) -
\begin{pmatrix}
0 & q_2(t) \\
q_3(t) & 0
\end{pmatrix}
y(t), \ y\in D(L_{dir}),
$$
where $\nu$ is defined by \eqref{dir6} and  $q_2$, $q_3$ -- by (\ref{dir7}).

The similarity transfrom $W_{dir}$ of $L_{dir}$ into $\widetilde{L}_{dir}$ is given by
\begin{equation}\label{dir10}
(W_{dir}y)(t)=
\begin{pmatrix}
e^{i\varphi(t)} & 0 \\
0 & e^{i\psi(t)}
\end{pmatrix}
y(t),
\end{equation}
$t\in[0, \omega]$, $y\in D(L_{dir})$, where $\varphi$, $\psi$ are defined by (\ref{dir5}).
\et

Observe that with this similarity transform the free operator satisfies $\widetilde L_{dir}^0={L}_{dir}^0- \nu I$. Therefore,  one can use the results in \cite{BDS11} to obtain all of the theorems for $L_{dir}$ stated in this paper, and we shall omit the corresponding proofs.

The following corollary is immediate from Theorem~\ref{dirth4}.
\bc\label{dirth5}
Assume that the potential $P$ is such that $p_2=p_3=0$. Then the operator $L_{dir}$ is similar to the operator $L^0_{dir} - \nu I$, where $\nu$ is defined by \eqref{dir6}. In particular, the eigenvalues of $L_{dir}$ are
$\lambda_n = \frac{\pi n}\omega -\nu$, $n\in\mathbb{Z}$, and the eigenfunctions $\widetilde s_n$ of $L_{dir}$ are the images of the eigenfunctions  of the free operator $L_{dir}^0$ under $W_{dir}$:
\[
\widetilde s_n = \frac1{\sqrt 2}
\begin{pmatrix}
e^{i\left(\varphi(t) - \frac{\pi n}{\omega}t\right)} \\
e^{i\left(\psi(t) + \frac{\pi n}{\omega}t\right)} 
\end{pmatrix},
 \ n\in\ZZ.
\]
\ec

\subsection{The preliminary similarity transforms for $L_{per}$ and $L_{ap}$.}

In this subsection, we exhibit new similarity transforms that are especially suited for $bc\in\{per, ap\}$.
Alternatively, one could use the similarity transform given by \eqref{dir10} with the corresponding domains for operators $L_{per}$ and $L_{ap}$. In the case when $\widehat{p}_1(0)$ and $\widehat{p}_4(0)$ are not both real, however, this would lead to a free operator that is not normal. The use of the method of similar operators would then become unnecessarily complicated.

The following theorem is the result of a straightforward computation, which we omit for the sake of brevity.

\bt\label{dirth6}
For $bc\in\{per, ap\}$, we have
\begin{equation}\label{dir11}
\widetilde{L}_{bc}=W_{bc}^{-1}L_{bc}W_{bc},
\end{equation}
where $W_{bc}\in B(L^2[0, \omega])$ and $\widetilde{L}_{bc}$ are given by
 \begin{equation}\label{dir13}
(W_{bc}y)(t)=
\begin{pmatrix}
e^{i(\varphi(t)-\frac{\theta}{\omega}t)} & 0 \\
0 & e^{i(\psi(t)-\frac{\theta}{\omega}t)}
\end{pmatrix}
y(t)  
\end{equation}
 and
\begin{equation}\label{dir12}
(\widetilde{L}_{bc}y)(t)=i
\begin{pmatrix}
1 & 0 \\
0 & -1
\end{pmatrix}
\frac{dy}{dt} -
\begin{pmatrix}
\widehat{p}_1(0) & 0 \\
0 & \widehat{p}_4(0)
\end{pmatrix}
y(t)-
\begin{pmatrix}
0 & q_2(t) \\
q_3(t) & 0
\end{pmatrix}
y(t),
\end{equation}
$y\in D(L_{bc})$, $t\in [0, \omega]$, with
\begin{equation}\label{dir14}
D(\widetilde{L}_{bc})=D(L_{bc}).
\end{equation}
In the above formulae, $\theta$ is given by (\ref{dir6}), $\varphi$ and $\psi$~-- by
(\ref{dir5}), and $q_2$ and $q_3$ -- by (\ref{dir7}).
\et

From Theorem~\ref{dirth6}, it follows that the study of the operators $L_{bc}$ is reduced to that of $\widetilde{L}_{bc}$, $bc\in\{per, ap\}$.
We write
$$
\widetilde{L}_{bc}=\widetilde{L}_{bc}^P-Q: D(\widetilde{L}_{bc})=D(L_{bc})\subset L_2[0, \omega]\to L_2[0, \omega],
$$
where the free (unperturbed) operator $\widetilde{L}_{bc}^P$ is given by
\begin{equation}\label{dir15}
(\widetilde{L}_{bc}^P y)(t)=i
\begin{pmatrix}
1 & 0 \\
0 & -1
\end{pmatrix}
\frac{dy}{dt} -
\begin{pmatrix}
\widehat{p}_1(0) & 0 \\
0 & \widehat{p}_4(0)
\end{pmatrix}
y(t), \quad y\in D(\widetilde{L}_{bc}), \quad t\in [0, \omega],
\end{equation}
and the perturbation $Q: D(L_{bc})\subset L_2[0, \omega]\to L_2[0, \omega]$ -- by
\begin{equation}\label{dir16}
(Qy)(t)=
\begin{pmatrix}
0 & q_2(t) \\
q_3(t) & 0
\end{pmatrix}
y(t), \quad y\in D(\widetilde{L}_{bc}), \quad t\in [0, \omega],
\end{equation}
with $q_2$ and $q_3$ defined by (\ref{dir7}).

The following corollary is immediate from Theorem~\ref{dirth6}.

\bc\label{dirth7}
Assume that the potential $P$ is such that $p_2=p_3=0$. Then the operator $L_{bc}$ is similar to the operator $\widetilde L_{bc}^P$ given by \eqref{dir15}, and we have
\begeq\label{spdec}
\s(L_{bc}) = \bigcup_{n\in\ZZ} \s_n, 
\eq
where
\ben
\item[(a)] in the case $bc=per$,   $\sigma_n=\left\{\frac{2\pi n}{\omega}-\widehat{p}_1(0), \frac{2\pi n}{\omega}-\widehat{p}_4(0)\right\}$, $n\in\mathbb{Z}$, and the corresponding eigenvectors of $L_{bc} = L_{per}$ are given by
\begin{equation}\label{eigper}
\bs
g_{n}^1(t) &=
\begin{pmatrix}
e^{-i\left((\frac{2\pi n}{\omega}-\widehat{p}_1(0))t+\int_0^t {p}_1(\tau)d\tau\right)} \\
0
\end{pmatrix}
\mbox{ and } \\
g_{n}^2(t) &=
\begin{pmatrix}
0 \\
e^{i\left((\frac{2\pi n}{\omega}-\widehat{p}_4(0))t+\int_0^t {p}_4(\tau)d\tau\right)}
\end{pmatrix},  \ n\in\mathbb{Z}, \ t\in [0, \omega];
\end{split}
\eq

\item[(b)] in the case $bc=ap$, $\sigma_n=\left\{\frac{\pi(2n+1)}{\omega}-\widehat{p}_1(0), \frac{\pi(2n+1)}{\omega}-\widehat{p}_4(0)\right\}$,
$n\in\mathbb{Z}$, and the corresponding eigenvectors of $L_{bc} = L_{ap}$  are given by
$$
\bar g_{n}^1(t)=  e^{-\frac{\pi i t}\omega}g_{n}^1(t)
\mbox{ and }
\bar g_{n}^2(t)= e^{\frac{\pi i t}\omega}g_{n}^2(t)
, \ n\in\mathbb{Z},\  t\in [0, \omega],
$$
where $g_n^1$ and $g_n^2$ are given by \eqref{eigper}.
\een
\ec

\section{Constructing the transforms $J$ and $\Gamma$ for Dirac operators.}\label{transext}

In this section, we adapt the construction from Section \ref{genJG} to obtain the transforms of the method of similar operators for 
$\widetilde L_{bc}$, $bc\in\{per, ap\}$, defined by \eqref{dir12}. We then show that inside the family of admissible triplets we get in the process, there will be at list one that satisfies Assumption \ref{baskpred1} for a given potential matrix that we consider in this paper. 

As in Section \ref{genJG}, we choose $\mathcal M = \mathfrak S_2(\HH)$ to be the space of admissible perturbations. Since we deal with Dirac operators, we have $\HH = L^2([0,\omega], \CC^2)$. For each potential matrix $P$, we 
 will exhibit a pair of families of transforms -- $\{J_k = J_k^P$, $k\in\ZZ_+\}$ and $\{\Gamma_k = \Gamma^P_k$, $k\in\ZZ_+\}$ --  that form admissible triplets together with $\mathcal M$.
Assuming  \eqref{badcond},
the  family  $\{J_k$, $k\in\ZZ_+\}$ will be independent of $P$ and we will  omit the superscript. The second family,  however, does depend on $P$ and we will only omit the superscript in cases when no confusion may arise. We shall write $\{\Gamma^0_k$, $k\in\ZZ_+\}$  for the transforms corresponding to $\widetilde L^0_{bc} = L^0_{bc}$.

\begin{rem}\label{rembad2}
If \eqref{badcond} does not hold, 
the construction from Section \ref{genJG} leads to
different families $\{J_k^P\}$ and $\{\Gamma_k^P\}$
for each integer value of $r = \frac{\omega\beta}{2\pi} = \frac\omega{2\pi}\left(\widehat{p}_4(0) - \widehat{p}_1(0)\right)$.
All of these cases, however, can be treated following the blueprint of \cite{BDS11}. As the corresponding constructions contain no new insights, we shall only state the results and omit the proofs. 
\erem

To use the construction from Section \ref{genJG}, we will need the spectral resolutions of the identity akin to the ones in Example \ref{direx}. As before, in the case $bc = per$, we let
$\l_n = \frac{2\pi n}{\omega}$, $n\in\mathbb{Z}$, and  in the case $bc = ap$, we let
$\l_n = \frac{\pi (2n+1)}{\omega}$, $n\in\mathbb{Z}$. The spectral components $\s_n = \s_n^P$ will then be as in Corollary \ref{dirth7}: 
$
\sigma_n=\left\{\l_n-\widehat{p}_1(0),\, \l_n-\widehat{p}_4(0)\right\}.
$
Notice that $\s_n$ is a two-element set if $\widehat{p}_1(0)\neq\widehat{p}_4(0)$ and a singleton otherwise. Moreover, these components form a (disjoint) partition of the spectrum of
$\widetilde L^P_{bc}$ because of our assumption \eqref{badcond} on $P$. 
Next,
we let $\{P_n, n\in\ZZ\}$ be the resolution of the identity that corresponds to the spectral decomposition \eqref{spdec} of $\widetilde L^P_{bc}$  with these $\s_n$. As in Example \ref{direx},
we shall also  use the spectral projections $P_{(k)} = \sum\limits_{|n|\le k} P_n$,  and the corresponding coarser resolutions of the identity given by $\{P_{(k)}\}\cup\{P_n, |n|> k\}$, $k\in \ZZ_+$. Notice that 
these projections are independent of the potential matrix $P$. In fact, each $P_n$, $n\in\ZZ$, is an orthogonal projection onto the eigenspace spanned by the vectors
\begeq\label{basis}
e_n^1=
\begin{pmatrix}
e_{-n} \\
0
\end{pmatrix}
\mbox{ and }
e_n^2 =
\begin{pmatrix}
0 \\
e_n
\end{pmatrix},
\quad e_n(t)
=e^{i\lambda_nt}, \quad t\in [0, \omega].
\eq

With the above notation, \eqref{tJ2} and \eqref{tG2} give us the transforms $J_0$ and $\Gamma_0^0$, respectively. The first equation in \eqref{Jk} then gives us the transforms $J_k$ and the first equation in
\eqref{Gk} -- the transforms $\Gamma_k^0$, $k\in\ZZ_+$. To define the transforms in the general case, we had to choose between two options. The first choice was to follow the abstract scheme outlined in Section \ref{genJG}. This, however,   made the family $\{J_k\}$ dependent on $P$ and the subsequent calculations became insurmountably difficult. Thus, we chose to
pursue the second option, where we kept the family $\{J_k\}$ intact and defined $\{\Gamma_k^P\}$ as a perturbation of $\Gamma_k^0$: we let $\Gamma_k^P = \Gamma_k^0+\Delta_k^P$. The formula for $\Delta_k^P$ arose from the key Property \ref{keyprop} in Definition \ref{baskdef8} of admissible triplet.

To exhibit the formula for $\Delta_k^P$,
it is convenient to use the following notation for the matrix elements of operators $X: D(X)\subseteq\HH\to \HH$. We will denote by $X_{mn}$ the matrix element 
that corresponds to the operator $P_{m}XP_n$. Clearly, we have
\begeq\label{Xmatr}
X_{mn} = 
\begin{pmatrix}
x_{mn}^{11} & x_{mn}^{12}\\
x_{mn}^{21} & x_{mn}^{22}
\end{pmatrix}
 = 
\begin{pmatrix}
\la Xe_{n}^1, e_m^1\ra & \la Xe_{n}^2, e_m^1\ra\\
\la Xe_{n}^1, e_m^2\ra & \la Xe_{n}^2, e_m^2\ra
\end{pmatrix},
\eq
where the basis elements $e_n^{1,2}$, $n\in\ZZ$, are given by \eqref{basis}. 
Observe that the matrix elements of the operator $Q$ given by \eqref{dir16} satisfy
\begeq\label{qmn}
Q_{mn} = 
\begin{pmatrix}
0 &  \widehat{q}_2(-m-n-\epsilon_{bc})\\
\widehat{q}_3(m+n+\epsilon_{bc}) &  0
\end{pmatrix},\quad m,n\in\ZZ,
\eq
where $q_k$, $k=2,3$, are given by \eqref{dir7}, $bc\in\{per, ap\}$,  $\epsilon_{per} = 0$, and $\epsilon_{ac} = 1$.

For $X \in \mathfrak S_2(\HH)$ given by \eqref{Xmatr}, we  let $\beta = \widehat{p}_1(0) -  \widehat{p}_4(0)$ as in \eqref{badcond} and
\begeq\label{deltap} 
(\Delta_0^P X)_{mn} = \frac{ \beta}{\l_m-\l_n}
\begin{pmatrix}
0 & \frac{ x_{mn}^{12}}{\l_m-\l_n-\beta}\\
\frac{-x_{mn}^{21}}{\l_m-\l_n+\beta} & 0
\end{pmatrix},\ m,n\in\ZZ, m\neq n.
 \eq
Naturally, we also define $(\Delta_0^P X)_{nn} = 0$, $n\in\ZZ$.
We note that $\Delta_0^P$ is well defined because of our technical assumption \eqref{badcond}. We also observe that 
since
$\l_m-\l_n = \frac{2\pi}\omega(m-n)$, we have $\Delta_0^P\mathfrak S_2(\HH)\subseteq \mathfrak S_2(\HH)$.
With the above definition, a straightforward computation shows that    Property \ref{keyprop} in Definition \ref{baskdef8} is, indeed, satisfied for $\Gamma_0^P = \Gamma_0^0+
\Delta_0^P$. It remains to let $\Gamma_k^P X = \Gamma_0^P X- P_{(k)}(\Gamma_0^P X)P_{(k)}$ as in \eqref{Gk}. Observe that the matrix elements of $\Gamma_k^PX$ satisfy
\begeq\label{GkP}
(\Gamma_k^P X)_{mn} =
\begin{cases}
\begin{pmatrix}
\frac{x_{mn}^{11}}{\l_m-\l_n} & \frac{x_{mn}^{12}}{\l_m-\l_n+\beta}\\
\frac{x_{mn}^{21}}{\l_m-\l_n-\beta} & \frac{x_{mn}^{22}}{\l_m-\l_n}
\end{pmatrix}, & \max\{|m|, |n|\}> k, m\neq n;\\
0, & \mbox{otherwise}. 
\end{cases}
\eq
We remark that the transforms $\Gamma_k^P$ do not come from the general scheme in Section \ref{genJG} if $\beta\neq 0$. The difference, however, is restricted to the main (block) diagonal in the block-matrix representation of $\Gamma_k^PX$. 

The advantage of the use of the transforms $J_k$ and $\Gamma_k^0$ in the definition of $\Gamma_k^P$ is in the fact that these transforms were studied in \cite{BDS11, BP17}. Thus, we can simply recall their properties.

\bl\label{bask100}\cite{BDS11}.
The  following properties hold.
\ben
\item\label{BHext} 
The transforms $J_k$ and $\Gamma_k^0$, $k\in \ZZ_+$,  are well defined for $X \in B(\HH)$ by \eqref{tJ2}, \eqref{tG2}, and the first equations in \eqref{Jk} and \eqref{Gk}.
\item The transforms $J_k$, $k\in \ZZ_+$,  are idempotents and $$\|J_k\|_{B(B(\HH))}
= \|J_k\|_{B(\mathfrak S_2(\HH))} = 1.$$
\item The transform $\Gamma_0^0$   satisfies 
$$\|\Gamma_0^0\|_{B(B(\HH))}\le\frac{\omega}{4} \mbox{ and }  
\|\Gamma_0^0\|_{B(\mathfrak S_2(\HH))} \le\frac\omega{2\pi}.$$
\item\label{adprop41}
For  $X\in B(\HH)$, we have $\Gamma_k^0 X \in D({\rm{ad}}_{L^0_{bc}})$ and
\[
{\rm{ad}}_{L^0_{bc}}(\Gamma_k^0 X) =  {L^0_{bc}}(\Gamma_k^0 X) - (\Gamma_k^0 X){L^0_{bc}} = X - J_k X.
\]
\een
\el
 
 From the above result we derive the following properties of  the transforms $\Gamma_k^P$, $k\in \ZZ_+$.

\bl\label{prop1GP} 
The  following properties hold.
\ben
\item\label{BHGext} 
The transforms   $\Gamma_k^P = \Gamma_k^0+\Delta_k^P$, $k\in \ZZ_+$,  are well defined for $X \in B(\HH)$ by \eqref{tG2},  \eqref{deltap}, and the first equation in \eqref{Gk}. 
\item\label{adprop22} The transforms $\Gamma_k^P$, $k\in \ZZ_+$,  belong to ${B(B(\HH))}$ as well as ${B(\mathfrak S_2(\HH))}$ and 
$$ 
\lim_{k\to\infty}\|\Gamma_k^P X\|_2=0 \mbox{ for all } X\in\mathfrak S_2(\HH).$$
\item\label{adprop11}
For  $X\in B(\HH)$, we have $\Gamma_k^P X \in D({\rm{ad}}_{\widetilde L^P_{bc}})$ and
\begeq\label{keyeq}
{\rm{ad}}_{\widetilde L^P_{bc}}(\Gamma_k^P X) =  {\widetilde L^P_{bc}}(\Gamma_k^P X) - (\Gamma_k^P X){\widetilde L^P_{bc}} = X - J_k X.
\eq
\een
\el
\bpf
In view of Lemma \ref{bask100}, to prove the first property, one needs to verify that \eqref{deltap} defines a bounded operator 
$\Delta_0^PX \in B(\HH)$ for any $X\in B(\HH)$.
First, recall from \cite{BDS11} that an arbitrary $X\in B(\HH)$ satisfies  $$Xx = \lim_{n\to\infty}\sum_{\ell= -n}^n\left(1-\frac{|\ell|}n\right) X_\ell x, \quad x\in\HH,$$
where $X_\ell$, $\ell\in\ZZ$, are the matrix diagonals of $X$.   
Secondly, for a given $\ell\neq 0$, and $x = \sum\limits_{n\in\ZZ}\left(\la x, e_n^1\ra e_n^1+\la x, e_n^2\ra e_n^2\right) = \sum\limits_{n\in\ZZ} (x_n^1 + x_n^2)\in \HH$, we have
 \[
 (\Delta_0^PX_\ell)x = \frac{\omega\beta}{2\pi} \sum_{n\in\ZZ} \frac1\ell
 \left(
\frac{x_{n+\ell,n}^{12}x_n^2}{\frac{2\pi\ell}{\omega}-\beta} - \frac{x_{n+\ell,n}^{21}x_n^1}{\frac{2\pi\ell}{\omega}+\beta} 
 \right)
 \]
from \eqref{deltap}. Since $\Delta_0^PX_0 = 0$, 
it follows that the sequence $$\{\|\Delta_0^PX_\ell\|_{B(\HH)}\}_{\ell\in\ZZ}$$ is summable, yielding Property \ref{BHGext} and $\|\Delta_0^P\|_{B(B(\HH))}<\infty$. The remaining part of Property \ref{adprop22} follows from \eqref{GkP}. More precisely, we have
\begeq\label{ng2}
\|\Gamma_k^P\|_{B(\mathfrak S_2(\HH))} =  
\max_{\substack{m, n\in\mathbb{Z}\\ m\neq n}} d_{mn}^{-1}
= :\delta^P, \quad k\in\mathbb{Z}_+,
\eq
where $d_{mn}=\mathrm{dist}_{m\ne n}\,(\sigma_m, \sigma_n)$, $m, n\in\mathbb{Z}$, is the distance between distinct spectral components $\sigma_m$ and $\sigma_n$ of the operator $\widetilde L_{bc}^P$. Observe that $d_{mn} = d_{nm}=d_{m-n,0}$, $m,n\in\ZZ$, and \eqref{GkP} implies
\[
\delta^P= \max \left\{\frac\omega{2\pi},\max_{\ell\in\ZZ\setminus\{0\}}  \frac{\omega}{\left|{2\pi}\ell-\beta\omega\right|}\right\}.
\]

Finally, the third property can be verified by direct computation; in fact, \eqref{deltap} was derived from \eqref{keyeq}.
\epf

From the above two lemmas and the preceding discussion we conclude that the following result holds.

\bc
The triplets $(\mathfrak S_2(\HH), J_k, \Gamma_k^P)$, $k\in\ZZ_+$, are admissible for $\widetilde L_{bc}^P$.
\ec

Next, to verify the properties in Assumption \ref{baskpred1}, we need to extend the transforms $J_k$ and $\Gamma_k^P$, $k\in\ZZ_+$, to $\mathfrak{L}_{\widetilde L^P_{bc}}(\mathcal{H})$. We recall that  $D(\widetilde L_{bc}^P) = D(L_{bc}) = D(L_{bc}^0)$, and that the spaces  $\mathfrak{L}_{\widetilde L^P_{bc}}(\mathcal{H})$ and $\mathfrak{L}_{ L^0_{bc}}(\mathcal{H})$ have equivalent norms. This allows us to define  the transforms so that $J_k$, $\Gamma^P_k$:  	
$\mathfrak L_{\widetilde L_{bc}^P}(\HH) \to \mathfrak L_{\widetilde L_{bc}^P}(\HH)$, $k \in \ZZ_+$, in the following way. Given $\l \in \rho({L_{bc}^0})$, we let
\begeq\label{Jext}
J_kX = J_k(X({L_{bc}^0}-\l I)^{-1}) ({L_{bc}^0}-\l I), \quad X\in  \mathfrak L_{\widetilde L_{bc}^P}(\HH),
\eq
and
\begeq\label{Gext}
\Gamma_k^PX = \Gamma_k^P(X({L_{bc}^0}-\l I)^{-1}) ({L_{bc}^0}-\l I), \quad X\in  \mathfrak L_{\widetilde L_{bc}^P}(\HH).
\eq
We observe that these extensions do not depend on the choice of $\l\in \rho({L_{bc}^0})$. 
Moreover, if $x \in D({L_{bc}^0})$, the formulas \eqref{tJ2},  \eqref{tG2},  \eqref{deltap}, and the first equations in \eqref{Jk} and  \eqref{Gk} may still be used for computation.

Yet again, we shall derive the properties of $J_k, \Gamma_k^P: \mathfrak{L}_{\widetilde L^P_{bc}}(\mathcal{H}) \to \mathfrak{L}_{\widetilde L^P_{bc}}(\mathcal{H})$ using analogous results in \cite{BDS11}.

\bl\label{dirlh3}\cite[Lemma 6]{BDS11}.
For the operator $Q$ given by \eqref{dir16}, we have  that the operators $J_k Q$, $\Gamma_k^0Q$, $Q\Gamma_k^0Q$, and $(\Gamma^0_k Q)J_kQ$, $k\in\mathbb{Z}_+$,
are Hilbert-Schmidt, i.e.~belong to $\mathfrak S_2(\HH)$.
\el

\bp\label{predass}
There exists $k\in\ZZ_+$ such that Assumption \ref{baskpred1} holds for the triplet
$(\mathfrak S_2(\HH), J_k, \Gamma_k^P)$ and the operator $Q$ given by \eqref{dir16}.
\ep

\bpf
Observe that it suffices to verify Properties \ref{propas2}, \ref{mprop}, and \ref{propas5} of Assumption \ref{baskpred1} for $\Gamma = \Delta_0^P$ and $B = Q$. Indeed, the remaining assertions will follow immediately from the definitions. Let us exhibit
the (non-zero) matrix elements of the relevant operators. To shorten the notation, we let
\begeq\label{quv}
u^{mn}_{bc} =  \widehat{q}_2(-m-n-\epsilon_{bc}) \mbox{ and } v^{mn}_{bc} =  \widehat{q}_3(m+n+\epsilon_{bc}),
\eq
where $q_k$, $k=2,3$, and $\epsilon_{bc}$ are as in \eqref{qmn}.
Using \eqref{deltap}, we have
\[
(\Delta_0^P Q)_{mn}  = \frac{\omega^2 \beta}{2\pi(m-n)}
\begin{pmatrix}
0 & \frac{ u^{mn}_{bc}}{2\pi(m-n)-\omega\beta}\\
\frac{-v^{mn}_{bc}}{2\pi(m-n)+\omega\beta} & 0
\end{pmatrix}, \ m\neq n\in\ZZ, 
\] 
\[
(Q\Delta_0^P Q)_{mn}  = \sum_{\substack{\ell\in\ZZ \\ {\ell\neq n}}}\frac{\omega^2 \beta}
{2\pi(\ell-n)}
\begin{pmatrix}
\frac{-u^{m\ell}_{bc}v^{\ell n}_{bc}}{2\pi(\ell-n)+\omega\beta} & 0\\ 
0 & \frac{ v^{m\ell}_{bc} u^{\ell n}_{bc}}{2\pi(\ell-n)-\omega\beta} 
\end{pmatrix}, \ \mbox{and}
\] 
\[
\bs
((&L_{bc}^0-\l I)(\Delta_0^PQ)(L_{bc}^0-\l I)^{-1})_{mn} = \\& 
\left[
\frac{\omega^2 \beta}
{2\pi n - \omega\l}+ \frac{\omega^2 \beta}
{2\pi(m-n)}\right] 
\begin{pmatrix}
0 & \frac{ u^{mn}_{bc}}{2\pi(m-n)-\omega\beta}\\
\frac{-v^{mn}_{bc}}{2\pi(m-n)+\omega\beta} & 0
\end{pmatrix}, \ m\neq n\in\ZZ,
\end{split}
\]
where $\l\in\rho(L_{bc}^0)$. Since $q_2, q_3\in L^2([0,\omega])$, we have $\widehat{q}_2$, $\widehat{q}_3\in\ell^2(\ZZ)$, and one can easily verify that all three of these operators are Hilbert-Schmidt via multiple applications of Cauchy-Schwarz inequality. For example, for $Q\Delta_0^P Q$, we have
\[
\begin{split}
\left\|Q\Delta_0^P Q\right\|_2^2 &= \sum_{m\in \ZZ}\sum_{\substack{n\in \ZZ\\ {n\neq m}}} \left\|(Q\Delta_0^P Q)_{mn} \right\|_2^2 \\&
=\sum_{m\in \ZZ}\sum_{\substack{n\in \ZZ\\ {n\neq m}}}\left| \sum_{\substack{\ell\in\ZZ \\ {\ell\neq n}}} \frac{\omega^2 \beta}
{2\pi(\ell-n)}\cdot\frac{u^{m\ell}_{bc}v^{\ell n}_{bc}}{2\pi(\ell-n)+\omega\beta}  \right|^2 \\&
+\sum_{m\in \ZZ}\sum_{\substack{n\in \ZZ\\ {n\neq m}}}\left| \sum_{\substack{\ell\in\ZZ \\ {\ell\neq n}}} \frac{\omega^2 \beta}
{2\pi(\ell-n)}\cdot\frac{u^{\ell n}_{bc}v^{m\ell}_{bc}}{2\pi(\ell-n)-\omega\beta}  \right|^2 =
\end{split}
\]
%
\[
\bs
 \left[\frac{\omega^2 \beta}{2\pi}\right]^2 &\sum_{m\in \ZZ}\sum_{\substack{n\in \ZZ\\ {n\neq m}}}\left(
\sum_{\substack{\ell\in\ZZ \\ {\ell\neq n}}} \sum_{\substack{\nu\in\ZZ \\ {\nu\neq n}}}
\frac{u^{m\ell}_{bc}v^{\ell n}_{bc}\overline{u^{m\nu}_{bc}v^{\nu n}_{bc}}}{(\ell-n)(\nu-n)(2\pi(\ell-n)+\omega\beta)(2\pi(\nu-n)+\omega\overline\beta)} \right.\\+& \left.
\sum_{\substack{\ell\in\ZZ \\ {\ell\neq n}}} \sum_{\substack{\nu\in\ZZ \\ {\nu\neq n}}}
\frac{u^{\ell n}_{bc}v^{m\ell}_{bc}\overline{u^{n\nu}_{bc}v^{m\nu}_{bc}}}{(\ell-n)(\nu-n)(2\pi(\ell-n)-\omega\beta)(2\pi(\nu-n)-\omega\overline\beta)}\right)
\le\\& 
\end{split}
\]
\[
\bs
&\frac{\omega^4 \beta^2}{4\pi^2} \left(\sum_{n\in \ZZ}
\sum_{\substack{\ell\in\ZZ \\ {\ell\neq n}}} \sum_{\substack{\nu\in\ZZ \\ {\nu\neq n}}}\sum_{\substack{m\in \ZZ\\ {m\neq n}}}
\frac{ |u^{m\ell}_{bc}|  |u^{m\nu}_{bc}| |v^{\ell n}_{bc}| |v^{\nu n}_{bc}|}{|\ell-n| |\nu-n| |2\pi(\ell-n)+\omega\beta| |2\pi(\nu-n)+\omega \beta|} 
\right.\\&+ \left.
\sum_{n\in \ZZ}\sum_{\substack{\ell\in\ZZ \\ {\ell\neq n}}} \sum_{\substack{\nu\in\ZZ \\ {\nu\neq n}}}\sum_{\substack{m\in \ZZ\\ {m\neq n}}}
\frac{|u^{\ell n}_{bc}||u^{n\nu}_{bc}| |v^{m\ell}_{bc}|  |v^{m\nu}_{bc}|}{|\ell-n| |\nu-n| |2\pi(\ell-n)-\omega\beta| |2\pi(\nu-n)-\omega\beta|}
\right)
\le 
\end{split}
\]
\[\bs
\frac{\omega^4 \beta^2}{4\pi^2} &\left(\|q_2\|_2^2\sum_{\substack{\ell\in\ZZ \\ {\ell\neq 0}}} \sum_{\substack{\nu\in\ZZ \\ {\nu\neq 0}}}\sum_{n\in \ZZ}
\frac{|v^{\ell,2n}_{bc}| |v^{\nu,2n}_{bc}|}{|\ell| |\nu| |2\pi\ell+\omega\beta| |2\pi\nu+\omega \beta|}\right.\\& \left.
+\|q_3\|_2^2\sum_{\substack{\ell\in\ZZ \\ {\ell\neq 0}}} \sum_{\substack{\nu\in\ZZ \\ {\nu\neq 0}}}\sum_{n\in \ZZ}
\frac{|u^{\ell,2n}_{bc}| |u^{\nu,2n}_{bc}|}{|\ell| |\nu| |2\pi\ell-\omega\beta| |2\pi\nu-\omega \beta|}
\right)
\le
\end{split}
\]
%
\[
\frac{\omega^4 \beta^2\|q_2\|^2_2\|q_3\|^2_2}{4\pi^2} \left(\sum_{\ell\in\ZZ\setminus\{0\}}\frac1{\ell^2}\right)
\left(\sum_{\nu\in\ZZ}\frac1{|2\pi\nu+\omega\beta|^2}+\sum_{\nu\in\ZZ}\frac1{|2\pi\nu-\omega\beta|^2}\right) <\infty.
\]
Properties \ref{propas2} and \ref{mprop} follow.

To prove Property \ref{propas5}, we write $Q(\widetilde L^P_{bc}-\l_\varepsilon I)^{-1} = MN_\varepsilon$, where 
$M: L^\infty \to L^2$ is given by $My = Qy$, $y\in L^\infty$, and $N_\varepsilon: L^2\to L^\infty$ -- by 
$N_\varepsilon x = (\widetilde L^P_{bc}-\l_\varepsilon I)^{-1}x$, $x\in L^2$. Clearly, $q_2, q_3\in L^2([0,\omega])$ implies that the operator $M$ is bounded. 
Now let us consider the operators $N_{i\ell}$, 
$\ell\in\mathbb{N}$. 
For $x=(x_1, x_2)\in\mathcal{H}$, we  have
$$
N_{i\ell} x=\left(\sum_{n\in\mathbb{Z}}\frac{\widehat{x}_1(n)e^{i\l_n t}}
{\l_n-\widehat p_{1}(0)-i\ell}, \sum_{n\in\mathbb{Z}}\frac{\widehat{x}_2(-n)e^{i\l_n t}}
{\l_n-\widehat p_{4}(0)-i\ell}\right).
$$
Therefore,
\begin{multline*}
\left\|N_{i\ell}x\right\|_\infty \leq \left(\sum_{n\in\mathbb{Z}}\frac{1}{\left(\l_n-
\mathrm{Re}\,\widehat{p}_1(0)\right)^2 + (\ell+\mathrm{Im}\,\widehat{p}_1(0))^2} + \right. \\ \left.+
\sum_{n\in\mathbb{Z}}\frac{1}{\left(\l_n-
\mathrm{Re}\,\widehat{p}_4(0)\right)^2 + (\ell +\mathrm{Im}\,\widehat{p}_4(0))^2}\right)^{1/2}\|x\|_2 \to 0
\end{multline*}
as $\ell\to\infty$, and the result follows.
\epf

We are now in position to apply Theorem \ref{baskth8}.

\bt\label{1trans}
There exists $k\in\mathbb{Z}_+$ such that  $\|\Gamma_k^PQ\|_2<1$, i.~e. the operator $I+\Gamma_k^PQ$ is invertible and $(I+\Gamma_k^PQ)^{-1}-I\in \mathfrak S_2(\HH)$. Moreover, the operator 
$\widetilde{L}_{bc}^P-Q$, $bc\in\{per, ap\}$, is similar to $\widetilde{L}_{bc}^P-{B}$,
${B}\in\mathfrak{S}_2(\mathcal{H})$, where
$$
{B}=J_kQ+(I+\Gamma_k^PQ)^{-1}(Q\Gamma_k^PQ-(\Gamma_k^PQ)J_kQ),
$$
and
$$
(\widetilde{L}_{bc}^P-Q)(I+\Gamma_k^PQ)=(I+\Gamma_k^PQ)(\widetilde{L}_{bc}^P- {B}).
$$
The operators $J_kQ$, $\Gamma_k^PQ$, $Q\Gamma_k^PQ$ are Hilbert-Schmidt and we have
that the operator ${B}\in\mathfrak{S}_2(\mathcal{H})$ satisfies
\begin{equation}\label{opB}
{B}=J_0Q+Q\Gamma_0^PQ
+C,
\end{equation}
where $C\in\mathfrak{S}_1(\mathcal{H})$.
\et
\begin{proof}
In view of Proposition \ref{predass} and Theorem \ref{baskth8}, we only need to verify the last assertion.
We have
\[
\begin{split}
C&=-(I+\Gamma_k^PQ)^{-1}(\Gamma_k^PQ)(Q\Gamma_k^PQ-(\Gamma_k^PQ)J_kQ) \\ &+Q\Gamma_k^PQ-Q\Gamma_0^PQ - (\Gamma_k^PQ)J_kQ 
+
J_kQ-J_0Q\\&=-(I+\Gamma_k^PQ)^{-1}(\Gamma_k^PQ)(Q\Gamma_k^PQ-(\Gamma_k^PQ)J_kQ) +C_2=C_1+C_2,
\end{split}
\]
where $C_2$ has finite rank  and $C_1\in\mathfrak{S}_1(\mathcal{H})$ as a product of two Hilbert-Schmidt operators.
\end{proof}


\section{Final similarity transform}\label{fsim}

To finish the proof of Theorem \ref{maint}, we would like to apply Theorem \ref{baskth6} to the operator $\widetilde L_{bc}^P - B$ constructed in Theorem \ref{1trans}. Unfortunately,  triplets $(\mathfrak S_2(\HH), J_m, \Gamma_m^P)$, $m\in \NN$, can be used only if $\|B\|_2$ is sufficiently small. Since we do not have any control over the norm of $B$, we would like to circumvent this issue.  
It is possible to do so if we choose a smaller space $\mathcal M$ of admissible perturbations. The space  $\mathcal M$  depends on the operator $B$ and its definition is rather involved. Roughly, it can be seen as a weighted Hilbert-Schmidt space.

 In what follows, it is assumed that the operator $B$ given by \eqref{opB} satisfies $B\neq  P_{(m)}BP_{(m)}$ for all $m\in\ZZ_+$. If the latter condition does not hold, the main results follow trivially from Theorem \ref{1trans}. 
To simplify the notation, in this section we shall write $A$ for $\widetilde L_{bc}^P$, $bc\in\{per, ap\}$.

Given $X\in\mathfrak S_2(\HH)$, we shall make use of the sequence $(\alpha_n(X))$  defined by
\begeq\label{alpha}
\alpha_n(X)=\|X\|_2^{-\frac{1}{2}}\max\left\{\left(\sum\limits_{\substack{|k|\ge  |n| \\
k\in\mathbb{Z}}}\|P_kX\|_2^2\right)^\frac{1}{4}, \left(\sum\limits_{\substack{|k|\ge  |n| \\
k\in\mathbb{Z}}}\|XP_k\|_2^2\right)^\frac{1}{4}\right\}, \quad n\in\mathbb{Z}.
\eq
It is easy to check that the above sequence has the following properties:
\ben
\item $\alpha_n(X)=\alpha_{-n}(X)$, $n\in\mathbb{Z}$.

\item $\lim\limits_{|n|\to\infty}\alpha_n(X)=0$, $n\in\mathbb{Z}$.

\item $\alpha_n(X)\leq  1$ for all $n\in\mathbb{Z}$.

\item $\alpha_n(X)\geq  \alpha_{n+1}(X)$, $n\geq  0$.

\item Given that $P_{(m)}XP_{(m)}\ne X$ for all $m\in\mathbb{Z}_+$, we have $\alpha_n(X)\ne 0$ for all $n\in\mathbb{Z}$, and
\begeq\label{weight}
\sum\limits_{n\in\mathbb{Z}}\frac{\|XP_n\|_2^2+\|P_nX\|_2^2}{(\alpha_n(X))^2} <\infty.
\eq
\een
We remark that the sequence $(\a_n(X))$ characterizes the decay of matrix elements of $X$ along the rows and columns of the matrix. In view of \eqref{weight}, one may conclude that any $X\in\mathfrak S_2(\HH)$ also belongs to a ``weighted Hilbert Schmidt space'' with a weight that depends on $X$. This is a manifestation of the fact that for any convergent series there is another one with a slower rate of convergence.  If $X = B$, we shall write simply $\a_n$ instead of $\a_n(B)$.

Next, we introduce the function $f = f_B: \sigma(A)\to\mathbb{R}_+$ given by
$$
f(\lambda)=\sum_{n\in\ZZ}\alpha_n\chi_{\s_n}(\l), \ \l\in\s(A),
$$
where $\s_n = \{\mu_n^1, \mu_n^2\}$ are, as before, the spectral components of the operator $A$, and $\chi_E$ denotes the characteristic function of a set $E$. In particular, if $bc = per$, we have $\mu_n^1 =  \frac{2\pi n}\omega -\widehat p_1(0)$ and $\mu_n^2 = \frac{2\pi n}\omega-\widehat p_4(0)$.
Using the functional calculus for unbounded normal operators \cite{DS88II}, we get that the operator
$$
f(A)=\sum_{n\in\mathbb{Z}}\alpha_nP_n
$$
belongs to $B(\mathcal{H})$ and $\|f(A)\|\leq  \max\limits_{n\in\mathbb{Z}}|\alpha_n|=1$.

We let $\M = \M(B)$ be the set of all operators  $X\in \mathfrak S_2(\HH)$ such that there 
exist operators $X_l, X_r\in\mathfrak S_2(\HH)$ satisfying
\begeq\label{xlxr}
X=X_lf(A) \quad\mbox{and}\quad X=f(A)X_r.
\eq
Observe that the matrix of the operator $f(A)$ is diagonal and the assumption $B\neq  P_{(m)}BP_{(m)}$ for all $m\in\ZZ_+$ implies that it  has no zeros on the main diagonal (see the last property of the sequence $\a_n$). Therefore, given $X \in\M$, the operators $X_l$ and $X_r$ are uniquely defined by \eqref{xlxr}. Moreover, $\M$ is a Banach space with the norm
$\|X\|_* = \max\{\|X_l\|_2, \|X_r\|_2\}$ and 
\begeq\label{Membed}
\|X\|_2 = \|X_lf(A)\|_2 = \|f(A)X_r\|_2 \le \|X\|_*.
\eq

From \eqref{weight}, we also deduce that $B\in\M(B)$ with 
$$B_l = \sum_{n\in\ZZ}\frac1{\a_n} BP_n \mbox{ and } B_r =  \sum_{n\in\ZZ}\frac1{\a_n} P_nB.$$

Next, we will show that for each $m\in\ZZ_+$ the transforms $J_m$ and $\Gamma_m^P$ defined in Section \ref{transext} form an admissible triplet for the operator $A=\widetilde L_{bc}^P$ together with $\M(B)$. We will also show that in this case the constants $\gamma = \g_m^P$ in Property \ref{adprop4} of Definition \ref{baskdef8} satisfy 
\begeq\label{limgam}
\lim\limits_{m\to \infty} \g_m^P = 0,
\eq
clearing the way for application of Theorem \ref{baskth6}.

The first and last properties of the admissible triplet in Definition \ref{baskdef8} follow from \eqref{Membed}.
Next, observe that
$$
J_m(X_lf(A))=(J_mX_l)f(A), \quad J_m(f(A)X_r)=f(A)(J_mX_r),
$$
$$
\Gamma_m^P(X_lf(A))=(\Gamma_m^PX_l)f(A), \quad \Gamma_m^P(f(A)X_r)=f(A)(\Gamma_m^PX_r),
$$
for $X_r$, $X_l\in\mathfrak{S}_2(\mathcal{H})$. It follows that the space $\M$ is invariant for $J_m$ and $\Gamma_m^P$ and, therefore,  the restrictions of the transforms to $\M$ are well defined. The second,  third, and fifth  properties of the admissible triplet follow. It remains to prove Property \ref{adprop4}. To do so, we shall make use of two more sequences:
$(\alpha_n')$ and $(\widetilde{\alpha}_n)$, $n\in\NN$. The first of them is defined by
$$
\alpha_{n+1}'=\max
\left\{\alpha_\ell{d_{j\ell}^{-1}}: \ell, j\in\mathbb{Z}, |\ell|\le n, |j|> n\right\}, \quad n\in\mathbb{Z}_+,
$$
where $d_{j\ell}=\mathrm{dist}
(\sigma_\ell, \sigma_j)$, $\ell, j\in\mathbb{Z}$, is, as before, the distance between the spectral components $\sigma_\ell$ and $\sigma_j$ of the operator $A$.
The second sequence is given by
\begeq\label{altil}
\widetilde{\alpha}_n=\delta^P\alpha_n+\alpha_n', \quad n\in\NN,
\eq
where $\delta^P$ is defined in \eqref{ng2}.
 Observe that  $\lim\limits_{n\to\infty}\alpha_n'= \lim\limits_{n\to\infty}\widetilde\alpha_n=0$.

An analog of the following lemma can be found in \cite{BDS11, BKR17}. We do, however, have to provide a new proof because the spectral components and the transforms are different in this paper. Moreover, the proof in this paper is more streamlined and provides better estimates.

\bl\label{dirlh7}
For any $m\in\mathbb{Z}_+$ and $X\in\mathfrak{S}_2(\mathcal{H})$ we have
$$
\max\left\{\|\Gamma_m^P(Xf(A))\|_2, \|\Gamma_m^P(f(A)X)\|_2\right\}\le  \widetilde{\alpha}_{m+1}\|X\|_2.
$$
\el
\begin{proof}
Let $Q_{(m)}=I-P_{(m)}$, $m\in\mathbb{Z}_+$. Then,
for $X\in\mathfrak{S}_2(\mathcal{H})$, we have 
$$
\Gamma_m^P(Xf(A))=\Gamma_m^P(Xf(A)Q_{(m)}) + \Gamma_m^P(Q_{(m)}Xf(A)P_{(m)}).
$$
Using \eqref{GkP} and the matrix representation \eqref{Xmatr} of $X$, we can write the matrix elements of the above two operators as
\[
(\Gamma_m^P (Xf(A)Q_{(m)}))_{j\ell} = 
\begin{cases}
\begin{pmatrix}
\frac{\a_\ell x_{j\ell}^{11}}{\l_j-\l_\ell} & \frac{\a_\ell x_{j\ell}^{12}}{\l_j-\l_\ell+\beta}\\
\frac{\a_\ell x_{j\ell}^{21}}{\l_j-\l_\ell-\beta} & \frac{\a_\ell x_{j\ell}^{22}}{\l_j-\l_\ell}
\end{pmatrix}, & |\ell|> m, \ell\neq j;\\
0, & \mbox{otherwise}; 
\end{cases}
\]
and
\[
(\Gamma_m^P(Q_{(m)}Zf(A)P_{(m)}))_{j\ell}= 
\begin{cases}
\begin{pmatrix}
\frac{\a_\ell x_{j\ell}^{11}}{\l_j-\l_\ell} & \frac{\a_\ell x_{j\ell}^{12}}{\l_j-\l_\ell+\beta}\\
\frac{\a_\ell x_{j\ell}^{21}}{\l_j-\l_\ell-\beta} & \frac{\a_\ell x_{j\ell}^{22}}{\l_j-\l_\ell}
\end{pmatrix}, & |\ell|\le m, |j|> m;\\
0, & \mbox{otherwise}. 
\end{cases}
\]
It follows that $$\|\Gamma_m^P (Xf(A)Q_{(m)})\|_2\leq \delta^P\alpha_{m+1}\|X\|_2$$ 
and
\[
\|\Gamma_m^P(Q_{(m)}Xf(A)P_{(m)})\|_2 
\le 
\alpha_{m+1}'\|X\|_2,
\]
yielding the first of the desired inequalities. The second one is obtained in a similar fashion.
\end{proof}

Immediately from the above lemma, we have the following result.

\bc
For any $m\in\mathbb{Z}_+$, we have $\|\Gamma_m^P\|_{B(\M)} \le \widetilde\a_{m+1}$.
\ec

The following lemma concludes the proof of Property \ref{adprop4} in Definition  \ref{baskdef8}. The proof is exactly the same as in \cite{BDS11} and is, therefore, omitted.

\bl\cite[Lemma 4]{BDS11}.
For any $m\ge 0$ and $X, Y \in \M$, we have 
\[
\max\left\{\|X\Gamma_m^P Y\|_*, \|(\Gamma_m^PX)Y\|_*\right\}\le \widetilde\a_{m+1}\|X\|_*\|Y\|_*.
\]
\el

We summarize the facts obtained after the definition of the space $\M = \M(B)$ in the following proposition.

\bp
For any $m\ge 0$ we have that $(\M(B), J_m, \Gamma_m^P)$ is an admissible triplet for the operator $A$ with $\|J_m\|=1$ and $\gamma_m = \widetilde\a_{m+1}$.
\ep

We are now in position to apply Theorem \ref{baskth6}. In conjunction with Lemma \ref{basklh1}, we obtain the following Theorem \ref{mainres11}, which is, in fact, stronger, than Theorem \ref{maint}. To simplify the exposition, we collect all the necessary notation and hypotheses in the following assumption.

\begin{ass}\label{mainass} Assume the following.
\ben
\item An operator $L_{bc} = L_{bc}^0 - P$ is given by \eqref{dir1}.
\item 
$W_{bc}$ is given by  \eqref{dir10} for $bc = dir$ and \eqref{dir13} for $bc\in\{per, ap\}$.
\item $\widetilde{L}_{bc}=W_{bc}^{-1}L_{bc}W_{bc} = \widetilde{L}_{bc}^P-Q$, where $Q$ is given by \eqref{dir16}.
\item $\beta = \widehat{p}_1(0) -  \widehat{p}_4(0)$, $r = \frac{\omega\beta}{2\pi}$ as in \eqref{badcond}.
\item If $bc = dir$ or $bc\in\{per, ap\}$ and 
$r \notin \ZZ\setminus\{0\}$, let
$\{P_n, n\in\ZZ\}$ be the resolution of the identity consisting of the spectral projections of the operator $L_{bc}^0$.
\item If $bc\in\{per, ap\}$ and 
$r \in \ZZ\setminus\{0\}$, let
$\{P_n, n\in\ZZ\}$ be the resolution of the identity consisting of the spectral projections of the operator $\widetilde L_{bc}^P$. 
\item Let $\{J_n, n\in\ZZ_+\}$ be the family of transforms defined by \eqref{tJ2} and \eqref{Jk}.
\item If $bc = dir$ or $bc\in\{per, ap\}$ and 
$r  \in \ZZ$, let $\{\Gamma_n, n\in\ZZ_+\}$ be the family of transforms defined by \eqref{tG2} and \eqref{Gk}.
\item If  $bc\in\{per, ap\}$ and 
$r\notin\ZZ $, let $\{\Gamma_n, n\in\ZZ_+\}$ be the family of transforms defined by  \eqref{GkP}.
\item For an operator $B\in\mathfrak S_2(\HH)$, let $\mathcal M =\mathcal M(B)$ be the Banach space defined earlier in this section, and $\widetilde \a_n = \widetilde \a_n(B)$, $n\in\NN$, be the sequence given by \eqref{altil}.
\een
\end{ass} 

\bt\label{mainres11}
With the notation and hypotheses in Assumption \ref{mainass},
 there exists $k\in\mathbb{Z}_+$ such that  $\|\Gamma_kQ\|_2<1$, and the operator 
${L}_{bc}$  is similar to $\widetilde{L}_{bc}^P-{B}$,
 where
$$
{B}=J_kQ+(I+\Gamma_kQ)^{-1}(Q\Gamma_kQ-(\Gamma_kQ)J_kQ)\in\mathfrak{S}_2(\mathcal{H}).
$$
Moreover, the operators $J_kQ$, $\Gamma_kQ$, and $Q\Gamma_kQ$ are Hilbert-Schmidt, and we have 
$
{B}=J_0Q+Q\Gamma_0Q
+C,
$
where $C\in\mathfrak{S}_1(\mathcal{H})$.

Furthermore, 
there exists $m\in\mathbb{Z}_+$ such that  $4\widetilde\a_{m+1}\|B\|_* < 1$, and the operator ${L}_{bc}$  is similar to $\widetilde{L}_{bc}^P-{V}$ with
$$
L_{bc}W_{bc}(I+U)=W_{bc}(I+U)(\widetilde{L}_{bc}^P-V).
$$
In the above formula:
\begin{itemize}
\item $U = \Gamma_kQ +\Gamma_m X_*+ (\Gamma_kQ)\Gamma_mX_* \in \mathfrak S_2(\HH)$;
\item $X_*\in \M(B)$ is the unique fixed point of the function $\Phi$ given by \eqref{bask13} with $J = J_m$ and $\Gamma = \Gamma_m$ in the ball \eqref{ball};
 \item $V = J_m X_* = 
 J_0Q+J_0(Q\Gamma_0Q)+\widetilde C\in 
\mathfrak{S}_2(\mathcal{H})$;
\item $\widetilde C\in\mathfrak S_1(\HH)$. 
 \end{itemize}
 Moreover, the subspaces $\mathcal{H}_{(m)} = P_{(m)}\HH$ and $\mathcal{H}_n = P_n\HH$, $|n|>m$, are invariant for the operator 
 $V$, the dimension of $\HH_{(m)}$ is at most $2m+2$, the dimensions of $\mathcal{H}_n$, $|n|>m$, are at most $2$, and  the operator $L_{bc}$ is a $W_{bc}(I+U)$-orthogonal direct sum
$$
L_{bc}=W_{bc}(I+U)\left(\left(\widetilde{L}_{bc}^P-V\right)_{(m)}\oplus\left(\bigoplus_{|n|>m}\left(\widetilde{L}_{bc}^P-V\right)_n\right)\right)(I+U)^{-1}W_{bc}^{-1}.
$$
\et

\bpf
The case $bc = dir$ follows from   \cite[Theorem 6]{BDS11}. The case $bc\in\{per, ap\}$ and $r \in \ZZ\setminus\{0\}$ does not follow directly from there but may be obtained in the same way.

The remaining case: $bc\in\{per, ap\}$ and $r  \notin \ZZ\setminus\{0\}$, follows from the discussion preceding the statement of the theorem.

The representation of $L_{bc}$ as a direct sum was not explicitly discussed in \cite{BDS11}. It follows from $V = J_mX_*$ and Lemma \ref{basklh1}. The formula
\begeq\label{vrep}
V =  J_0Q+J_0(Q\Gamma_0Q)+\widetilde C
\eq
is obtained the same way as in \cite{BDS11} using the fact that $J_0\mathfrak S_1(\HH)\subseteq \mathfrak S_1(\HH)$.
\epf

\section{Asymptotic estimates of the spectrum}\label{sec7}

In this section, we operate under Assumption \ref{mainass} and use
the direct sum representation in Theorem \ref{mainres11}  to estimate the spectrum $\s(L_{bc}) = \s(\widetilde{L}_{bc}^P-V)$. We let $\s_{(m)} = \s\left((\widetilde{L}_{bc}^P-V)_{(m)}\right)$ and $\s_{n} = \s\left((\widetilde{L}_{bc}^P-V)_{n}\right)$, $|n|\ge m$.

First, observe that the operator $\widetilde{L}_{bc}^P$ has compact resolvent. Therefore, so do the operators $\widetilde{L}_{bc}^P-V$ and $L_{bc}$. Repeating the argument in \cite[Remark 2]{BDS11}, we get
\begeq\label{specrep}
\s(L_{bc}) = \s_{(m)}\cup\left(\bigcup_{|n|\ge m} \s_n\right).
\eq
We remark that an analog of the above result was proved in \cite{DM10}.

Secondly, we observe that the spaces $\HH_n$, $n\in\ZZ$, have dimensions at most $2$ and  the space $\HH_{(m)}$  -- at most $2m+2$. Therefore, the set $\s_{(m)}$ has at most $2m+2$ elements and each $\s_{n}$ -- at most $2$. 

In the following result, we summarize coarse estimates of $\s(L_{bc})$ that follow immediately from Theorem \ref{mainres11},  Lemma \ref{basklh1}, and \eqref{specrep}.

\bt\label{mainres1}
The following equations hold.
\begin{itemize}
\item $\s(L_{dir}) = \left\{\frac{\pi n}\omega-\nu+b^1_n: n\in\ZZ\right\}$, 
\item $\s(L_{per}) = \left\{\frac{2\pi n}\omega - \widehat p_1(0)+b^2_n, \frac{2\pi n}\omega - \widehat p_4(0)+b^3_n: n\in\ZZ\right\}$,
\item $\s(L_{ap}) = \left\{\frac{\pi(2 n+1)}\omega - \widehat p_1(0)+b^4_n, \frac{\pi(2 n+1)}\omega - \widehat p_4(0)+b^5_n: n\in\ZZ\right\}$;
\end{itemize}
where $\nu$ is given by \eqref{dir6} and $b^k = \{b^k_n: n\in\ZZ\} \in\ell^2(\ZZ)$, $k = 1,\ldots, 5$.
\et

To provide more accurate estimates of  the spectra, 
we need   to consider various cases in greater detail.

\subsection{The spectrum of $L_{dir}$.}

This case is covered by   \cite[Theorem 7]{BDS11}. The result is as follows.

\bt There exists $m\in\NN$ such that \eqref{specrep} holds. Moreover, 
we have $\s_{n} = \{\widetilde\l_n\}$, $|n| > m$, where
\begeq\label{ltdir}
\widetilde\l_n = \frac{\pi n}\omega - \nu- \theta_{2n}-\sum_{\ell\in\ZZ\setminus\{0\}}\frac{\theta^2_{\ell+2n}}\ell +c_n,
\eq
\[
\theta_n = \begin{cases}
\frac12\left(\widehat q_2(-\frac n2)+\widehat q_3(\frac n2)\right), & n\in 2\ZZ;\\
\frac12\left(\widehat{\eta}_2(-\frac{ n+1}2)+\widehat \eta_3(\frac {n+1}2)\right), & n\in2\ZZ+1; 
\end{cases},
\]
$\nu$ is given by \eqref{dir6}, $q_k$, $k=2,3$, -- by \eqref{dir7}, 
$\eta_2(t) = q_2(t)e^{-i\frac{\pi}\omega t}$, $\eta_3(t) = q_3(t)e^{i\frac{\pi}\omega t}$, and $c = (c_n)\in\ell^1(\ZZ)$.
\et

\subsection{The spectrum of $L_{bc}$, $bc\in\{per, ap\}$, with $r = \frac{\omega\beta}{2\pi} \notin \ZZ$.}
 In this case, we will use the following lemma, which is implied by the calculations in Example \ref{twotwo}. 

\bl\label{newlemma}
Consider a family of $2\times 2$ matrices $Z_n$, $n\in\NN$, given by
\[
Z_n = \left(
\begin{array}{cc}
z_n^{1} & 0 \\
0 & z_n^{2} 
\end{array}
\right)
-
 \left(
\begin{array}{cc}
b_n^{1} & b_n^{3} \\
b_n^{4} & b_n^{2} 
\end{array}
\right)
-
 \left(
\begin{array}{cc}
c_n^{1} & c_n^{3} \\
c_n^{4} & c_n^{2} 
\end{array}
\right),
\]
where $\varepsilon = \inf\limits_{n\in\ZZ}|z_n^1- z_n^2|> 0$,   the sequences $\{b_n^k: n\in\ZZ\}$, $k=1,2,3,4$, belong to $\ell^2(\ZZ)$, and the sequences 
$\{c_n^k: n\in\ZZ\}$, $k=1,2,3,4$, belong to $\ell^1(\ZZ)$. Then for all  $n\in\ZZ$ that are sufficiently large in absolute value, the matrix $Z_n$ is similar to
\[
\widetilde Z_n =  \left(
\begin{array}{cc}
z_n^{1}-b_n^1-d_n^1 & 0 \\
0 & z_n^{2}-b_n^2-d_n^2 
\end{array}
\right),
\]
where  the sequences 
$\{d_n^k: n\in\ZZ\}$, $k=1,2$, belong to $\ell^1(\ZZ)$.
\el

\bpf
Observe that for $n\in\ZZ$  that are sufficiently large in absolute value we have 
\[
\left\| \left(
\begin{array}{cc}
b_n^{1} & b_n^{3} \\
b_n^{4} & b_n^{2} 
\end{array}
\right)
+
 \left(
\begin{array}{cc}
c_n^{1} & c_n^{3} \\
c_n^{4} & c_n^{2} 
\end{array}
\right) \right\|_2 < \frac\varepsilon4.
\]
According to  Example \ref{twotwo}, we have that for such $n\in\ZZ$ the operator $Z_n$ is similar to
\[
\widetilde Z_n =  \left(
\begin{array}{cc}
z_n^{1}-b_n^1- c_n^1 -\nu_n^1 & 0 \\
0 & z_n^{2}-b_n^2- c_n^2 -\nu_n^2 
\end{array}
\right),
\]
where $\{\nu_n^k\}$, $k=1,2$, can be estimated using \eqref{bgx}.  In particular, they belong to $\ell^1$ as pointwise products of two $\ell^2$ sequences.
\epf

In view of the representation \eqref{vrep}, 
we may apply the above result to the sequence $\left(\widetilde{L}_{bc}^P-V\right)_n $, $|n| >m$. 

\bt\label{mainres1000}
There exists $m\in\NN$ such that \eqref{specrep} holds with
\[
\bs
 \s_n =\Bigg\{&\frac{\pi(2 n+\epsilon_{bc})}\omega -\widehat p_1(0) - \sum_{\ell\neq 2n}\frac{\omega\widehat q_2(-\ell-\epsilon_{bc})\widehat q_3(\ell+\epsilon_{bc})}{2\pi(\ell-2n)-\omega\beta}-d_n^1,  \\
&\frac{\pi(2 n+\epsilon_{bc})}\omega-\widehat p_4(0) - \sum_{\ell\neq 2n}\frac{\omega\widehat q_2(-\ell-\epsilon_{bc})\widehat q_3(\ell+\epsilon_{bc})}{2\pi(\ell-2n)+\omega\beta}-d_n^2 \Bigg\},\ |n|> m,
\end{split}
\]
where $q_k$, $k=2,3$, are given by \eqref{dir7}, $\epsilon_{bc}$ is as in \eqref{qmn}, and the sequences 
$\{d_n^k: n\in\ZZ\}$, $k=1,2$, belong to $\ell^1(\ZZ)$.
\et

\bpf
Since $r\notin\ZZ$, the transform $\Gamma_0 = \Gamma_0^P$ is given by \eqref{GkP}. Using \eqref{qmn} and
\begeq\label{qgqmn}
(Q\Gamma_0 Q)_{nn} =  \left(
\begin{array}{cc}
\sum\limits_{\ell\neq 2n}\frac{\omega\widehat q_2(-\ell-\epsilon_{bc})\widehat q_3(\ell+\epsilon_{bc})}{2\pi(\ell-2n)-\omega\beta} & 0 \\
0 & \sum\limits_{\ell\neq 2n}\frac{\omega\widehat q_2(-\ell-\epsilon_{bc})\widehat q_3(\ell+\epsilon_{bc})}{2\pi(\ell-2n)+\omega\beta} 
\end{array}
\right),
\eq
the result follows from \eqref{vrep} and  Lemma \ref{newlemma} via a direct computation. 
\epf

\subsection{The spectrum of $L_{bc}$, $bc\in\{per, ap\}$, with $r =\frac{\omega\beta}{2\pi} \in \ZZ$.}

Observe that in this case the matrix elements of $Q$ and $Q\Gamma_0 Q$ satisfy
\begeq\label{qmnr}
\bs
Q_{mn} &= 
\begin{pmatrix}
0 & \widehat q_2(-m-n-r-\epsilon_{bc}) \\
\widehat q_3(m+n+r+\epsilon_{bc}) & 0
\end{pmatrix}
\mbox{ and }
\\
(Q\Gamma_0 Q)_{nn} &=  \sum\limits_{\ell\neq 2n+r}\frac{\omega\widehat q_2(-\ell-\epsilon_{bc})\widehat q_3(\ell+\epsilon_{bc})}{2\pi(\ell-2n-r)}\left(
\begin{array}{cc}
1 & 0 \\
0 & 1 
\end{array}
\right),
\end{split}
\eq
where $Q$ is given by \eqref{qmn} and $\Gamma_0$ by \eqref{tG2}.

We will need the following lemma.

\bl\label{lemspec}
Let $Z_n = B_n +C_n$ be a sequence of matrices such that
\[
B_n =  \left(
\begin{array}{cc}
0 & b_n^{2} \\
b_n^{3} & 0 
\end{array}
\right) \mbox{ and } 
C_n =  \left(
\begin{array}{cc}
c_n^{1} & c_n^{3} \\
c_n^{4} & c_n^{2} 
\end{array}
\right),
\]
where  the sequences $\{b_n^k: n\in\ZZ\}$, $k=2,3$, belong to $\ell^2(\ZZ)$, and the sequences 
$\{c_n^k: n\in\ZZ\}$, $k=1,2,3,4$, belong to $\ell^1(\ZZ)$. Then
\begeq\label{srep}
\s(Z_n) = \left\{\sqrt{b_n^2b_n^3} + d_n^1, -\sqrt{b_n^2b_n^3} + d_n^2\right\},
\eq
where the sequences 
$\{d_n^k: n\in\ZZ\}$, $k=1,2$, belong to $\ell^{4/3}(\ZZ)$. Moreover,
if there exist $c, C > 0$ such that $cb_n^2 \le b_n^3 \le Cb_n^2$ for all $n\in\ZZ$ that are sufficiently large in absolute value, then 
the sequences 
$\{d_n^k: n\in\ZZ\}$, $k=1,2$, belong to $\ell^{1}(\ZZ)$.
\el

\bpf
Let $\s(Z_n) = \{\mu_n^+, \mu_n^-\}$. 
Observe that 
\begeq\label{trace}
\left\{\mu_n^+ +\mu_n^- = {\rm{tr}}(Z_n) = c_n^1+c_n^2: n\in\ZZ\right\}\in\ell^1. 
\eq
Hence, 
\[
 {\rm{det}}(Z_n) = \mu_n^+\mu_n^- = \mu_n^+(c_n^1+c_n^2)-(\mu_n^+)^2  =  \mu_n^-(c_n^1+c_n^2)-(\mu_n^-)^2  = - b_n^2b_n^3+\nu_n,
\]
where $\{\nu_n = c_n^1c_n^2-b_n^3c_n^3-b_n^2c_n^4 - c_n^3c_n^4: n\in\ZZ\}\in \ell^{2/3}$. Since $\{\mu_n^\pm\}\in\ell^2$,
it follows that $\{ \mu_n^\pm(c_n^1+c_n^2)\}$ and, consequently, $\{(\mu_n^\pm)^2-b_n^2b_n^3\}$  belong to $\ell^{2/3}$. 
Since
\[
\bs
&\left[\min\left\{\left|\mu_n^\pm +\sqrt{b_n^2b_n^3}\right|, \left|\mu_n^\pm-\sqrt{b_n^2b_n^3}\right|\right\}\right]^2 \le \\
&  
\left|\mu_n^\pm+\sqrt{b_n^2b_n^3}\right|\left|\mu_n^\pm-\sqrt{b_n^2b_n^3}\right|
= \left|(\mu_n^\pm)^2-b_n^2b_n^3\right|,
\end{split}
\]
we have
\[
\left\{\min\left\{\left|\mu_n^\pm  +\sqrt{b_n^2b_n^3}\right|, \left|\mu_n^\pm-\sqrt{b_n^2b_n^3}\right|\right\}: n\in\ZZ\right\}\in \ell^{4/3}.
\]
Using \eqref{trace} once again, we deduce \eqref{srep} with $\{d_n^k\}\in \ell^{4/3}$.

Let us now assume that $c, C > 0$ are such that $cb_n^2 \le b_n^3 \le Cb_n^2$ for all $n\in\ZZ$ that are sufficiently large in absolute value. Without loss of generality, we may also assume
that $b_n^2 \neq 0$ and the double inequalities hold for all $n\in\ZZ$. Then $r_n = \sqrt{b_n^3/b_n^2}$ is well defined and 
\begeq\label{rnest}
0 < \sqrt c \le r_n \le \sqrt C <\infty.
\eq
Therefore, the matrices
\[
U_{n} = 
\begin{pmatrix}
1 &  1\\
-r_n &  r_n
\end{pmatrix}
\quad \mbox{and}\quad
U_{n}^{-1} = \frac1{2r_n}
\begin{pmatrix}
r_n &  -1\\
r_n &  1
\end{pmatrix}
\]
are also well defined. Then the matrix
\[
\widetilde Z_N = U_n^{-1}Z_n U_n=  
\left(
\begin{array}{cc}
-\sqrt{b_n^2b_n^3} &  0\\
0 &  \sqrt{b_n^2b_n^3}
\end{array}
\right)  +U_n^{-1}\left(
\begin{array}{cc}
c_n^{1} & c_n^{2} \\
c_n^{3} & c_n^{4} 
\end{array}
\right) 
U_n
\]
has the same spectrum as $Z_n$, and the desired result follows from \eqref{rnest} and an argument similar to the proof of Lemma \ref{newlemma}. Indeed, the similarity transform
of Example \ref{twotwo} applies for the subsequence indexed by $n\in\ZZ$ such that
\begeq\label{morm}
\left\|U_n^{-1}\left(
\begin{array}{cc}
c_n^{1} & c_n^{2} \\
c_n^{3} & c_n^{4} 
\end{array}
\right) 
U_n\right\|_2 < \frac12\sqrt{|b_n^2b_n^3|}.
\eq
For the complementary subsequence where \eqref{morm} is not satisfied, we have that \eqref{srep} holds automatically because the sequence of norms in the left hand side of \eqref{morm} is in $\ell^1$ due to \eqref{rnest}.
\epf

The condition in the above lemma motivates the following definition.

\bd
For a given $r\in\ZZ$, the potential function $P$ given by \eqref{dir2} is called \emph{$r_{bc}$-balanced} if there exist $c, C > 0$ and $N\in\mathbb N$ such that 
\begeq\label{cC}
cq_2(-2n-r-\epsilon_{bc}) \le q_3(2n+r+\epsilon_{bc}) \le Cq_2(-2n-r-\epsilon_{bc})
\eq
for all $n\in\ZZ$ satisfying $|n|\ge N$. In the above double inequality $q_2$ and $q_3$ are given by \eqref{dir7}, and  $\epsilon_{bc}$ is as in \eqref{qmn}.
\ed

With the above definition, the following result is essentially immediate.

\bt\label{mainres10100}
There exists $m\in\NN$ such that \eqref{specrep} holds with
\begeq\label{snbad}
\bs
 \s_n =\Bigg\{&\frac{\pi(2 n+\epsilon_{bc})}\omega -\widehat p_1(0) -
\sum_{\ell\neq 2n+r}\frac{\omega\widehat q_2(-\ell-\epsilon_{bc})\widehat q_3(\ell+\epsilon_{bc})}{2\pi(\ell-2n-r)}
 \ -  \\ 
&\sqrt{\widehat q_2(-2n-r-\epsilon_{bc})\widehat q_3(2n+r+\epsilon_{bc})}-d_n^1,  \\
& \frac{\pi(2 n+\epsilon_{bc})}\omega -\widehat p_1(0) -
\sum_{\ell\neq 2n+r}\frac{\omega \widehat q_2(-\ell-\epsilon_{bc})
\widehat q_3(\ell+\epsilon_{bc})}{2\pi(\ell-2n-r)}
 \  +  \\ 
& \sqrt{\widehat q_2(-2n-r-\epsilon_{bc})\widehat q_3(2n+r+\epsilon_{bc})} 
-d_n^2 \Bigg\},\ |n|> m,
\end{split}
\eq
where the sequences 
$\{d_n^k: n\in\ZZ\}$, $k=1,2$, belong to $\ell^{4/3}(\ZZ)$. Moreover, if the potential function $P$ in \eqref{dir2} is $r_{bc}$-balanced then the sequences 
$\{d_n^k: n\in\ZZ\}$, $k=1,2$, belong to $\ell^{1}(\ZZ)$.
\et

\bpf
The result follows from \eqref{vrep}, \eqref{qmnr}, and Lemma \ref{lemspec}.
\epf

\brem
If $p_1(t)=p_4(t)\equiv 0$,  Theorem \ref{mainres10100} yields \cite[Theorem 7]{BDS11} in the case $bc \in\{per, ap\}$. Then $r = 0$ and the sequences $\{d_n^k\}$, $k=1,2$, are guaranteed to be in $\ell^{4/3}(\ZZ)$. The case of balanced potentials was not considered in \cite{BDS11}. A different estimate for the eigenvalues in the case of balanced potentials can be obtained in terms of the constants $c$ and $C$ in \eqref{cC}.
We cite \cite{BP17} for a similar approach used for a different kind of operator.
\erem

\section{Equiconvergence of spectral decompositions}\label{sec8}

In this section, we explore the consequences of Theorem \ref{mainres11} for equiconvergence of spectral decompositions. We continue to operate under Assumption \ref{mainass} and use the notation of Theorem \ref{mainres11}. We cite \cite[and references therein]{DM10, DM12JFA, DM12, LM14, LM15, MO12, S16t, S16, SS15, SS14} for preexisting results on spectrality and equiconvergence for Dirac operators. A few special cases of the results in this section appear in \cite{BDS11}. It turns out that the proofs in the general case are either very straightforward or very similar; we omit them for the sake of brevity.

\bd\cite[Definition 3]{BDS11}.
Assume that the spectrum $\s(A)$ of a linear operator $A: D(A)\subseteq \HH\to\HH$ satisfies
\begeq\label{decspec}
\s(A) =\bigcup_{n\in\ZZ} \s_k,
\eq
where $\s_k$ are mutually disjoint compact sets. Let $P_k$ be the Riesz projections corresponding to $\s_k$. The operator $A$ is said to be \emph{spectral with respect to the decomposition} \eqref{decspec} (or to possess \emph{generalized spectrality}) if $\{P_k: k\in\ZZ\}$ forms a resolution of the identity in $B(\HH)$. 
\ed

Immediately, from Theorems \ref{mainres11}, \ref{mainres1} and Lemma \ref{basklh1},
we get the following result (see also \cite{BDS11, DM10}).

\bt\label{mainresspec}
The operator $L_{bc}$ is spectral with respect to  decomposition \eqref{specrep}.
\et

Next, similarly to \cite{BDS11}, we compare the convergence of spectral decompositions for operators $L_{bc}$,  $\widetilde L^P_{bc}-Q$, and $\widetilde L^P_{bc} -V$. To this end,
we consider three families of Riesz projections. The first, $\{P_{(m)}\}\cup\{P_j: |j|> m\}$, 
corresponds to the spectral decompositions  of both $\widetilde L^P_{bc}$ and $\widetilde L^P_{bc} -V$. 
The second family of Riesz projections,  $\{\widetilde P_{(m)}\}\cup\{\widetilde P_j: |j|>m\}$, corresponds to the operator $\widetilde L^P_{bc}-Q$. From Theorem \ref{mainres11} and Lemma \ref{basklh1}, we have
\begeq\label{simproj}
\widetilde P_j = (I+U)P_j(I+U)^{-1},\ |j|>m.
\eq
The third family of Riesz projections,  $\{\bar P_{(m)}\}\cup\{\bar P_j: |j|>m\}$, corresponds to the operator $L_{bc}$. Again, from Lemma \ref{basklh1}, we have
\begeq\label{simproj2}
\bar P_j = W_{bc} \widetilde P_jW_{bc}^{-1},\ |j|>m.
\eq

Given a set $\Omega\subseteq \{j\in\ZZ: |j|> m\}$, we define 
\[
P(\Omega) = \sum_{j\in\Omega} P_j, \ \widetilde P(\Omega) = \sum_{j\in\Omega} \widetilde P_j, \mbox{ and }  \bar P(\Omega) = \sum_{j\in\Omega} \bar P_j.
\]
We also let $ \a(\Omega,X) = \max_{n\in\Omega} \a_n(X)$, where $\a_n(X)$ is defined by \eqref{alpha}.
The following lemma is proved the same way as \cite[Lemma 8]{BDS11}. 
\bl\label{BDSle8}
Given $X\in \mathfrak S_2(\HH)$ and $\Omega\subseteq \{j\in\ZZ: |j|> m\}$, we have
\[
\max\{\|P(\Omega)X\|_2, \|XP(\Omega)_2\} \le C(X) \a(\Omega,X)
\]
for some constant $C(X)$ that depends only on $X$.
\el

Immediately from \eqref{alpha}, we also have the following auxiliary result (see also \cite[Lemma 9]{BDS11}.

\bl\label{BDSle9}
Assume that $ 0\neq X\in \mathfrak S_2(\HH)$ and $\Omega\subseteq \{j\in\ZZ: |j|> m\}$.
\ben
\item
If $X = \sum_{j\in\ZZ} X_j$, $X_j\in \mathfrak S_2(\HH)$, and the series converges absolutely, then
\[
\|X\|_2^{\frac12} \a(\Omega,X) \le \sum_{j\in\ZZ} \|X_j\|_2^{\frac12} \a(\Omega,X_j).
\]
\item If $X = X_1\cdot \ldots \cdot X_\ell$, $X_j\in \mathfrak S_2(\HH)$, $1\le j\le \ell$, then
\[
\|X\|_2^{\frac12} \a(\Omega,X) \le \min_{1\le j\le \ell}\a(\Omega,X_j) \prod_{j=1}^\ell \|X_j\|_2^{\frac12} .
\]
\een
\el


The key result of this section is the following theorem (once again, see Assumption \ref{mainass} and Theorem \ref{mainres11} for the notation).

\bt
Given  $\Omega\subseteq \{j\in\ZZ: |j|> m\}$, we have
\[
\|\widetilde P(\Omega) - P(\Omega)\|_2 \le C(\a(\Omega, \Gamma_0 Q) +
\a(\Omega, \Gamma_0X_*)),
\]
where $C >0$ does not depend on $\Omega$.
\et

\bpf
The proof of \cite[Theorem 9]{BDS11} can be used as a blueprint.
\epf

Immediately from the above theorem and \eqref{simproj2} we get the following result about equiconvergence.

\bt\label{mainresequi}
We have
\[
\lim_{\ell\to\infty}\left\|\bar{{P}}_{(m)}+\sum_{m<|n|\le\ell}\bar{{P}}_n-W_{bc}^{-1}\left({P}_{(m)}+
\sum_{m<|n|\le\ell} P_n\right)W_{bc}\right\|_{2}=0.
\]
\et

\section{$C_0$-group generated by the Dirac operator}\label{groupsec}
We conclude this paper with a description of the $C_0$-group generated by the operator $iL_{bc}$ under Assumption \ref{mainass}. 
The following result is immediate from Theorems \ref{baskth7}, \ref{mainres11} and \ref{mainres1}, and Lemmas \ref{basklh1} and \ref{baskuskova_lh6}.
In its formulation, we continue to use the notation from Assumption \ref{mainass} and  Theorem \ref{mainres11}.

\bt\label{mainres3}
The operator $iL_{bc}$ generates a $C_0$-group of operators $T: \RR\to B(\HH)$ such that $T(t) = W_{bc}(I+U)\widetilde T(t) (I+U)^{-1}W_{bc}^{-1}$, where
\begeq\label{tildet}
\widetilde T(t) =e^{it\left(\widetilde L^P_{bc}-V\right)_{(m)}}\oplus\left(\bigoplus_{|n|> m} e^{it\left(\widetilde L^P_{bc}-V\right)_{n}}\right), \quad t\in\mathbb{R}.
\eq
\et

\bpf
Theorem \ref{mainres11} provides a $W_{bc}(I+U)$-orthogonal decomposition of the operator $L_{bc}$ and an orthogonal decomposition of the operator $\widetilde L^P_{bc}-V$.
From Theorem \ref{mainres1} we deduce that \eqref{bask_16'} is satisfied for $\widetilde L^P_{bc}-V$, and Lemma \ref{baskuskova_lh6} yields the group $\widetilde T$. It remains to cite Lemma \ref{basklh1}(4) to complete the proof. 
\epf

As a corollary, we obtain the following result, which can be used to estimate the rate of convergence of the
generalized Fourier series of  mild solutions of differential equations involving Dirac operators. We 
 let $Z = W_{bc}(I+U)$.

\bc\label{Fouriercon1}
For   any $\varphi \in \HH$ and  $n \in \NN$ with $n > m$, we have
\[
\left\|
T(t)\varphi - Z\widetilde T(t)P_{(n)}Z^{-1}\varphi
\right\| \le \|Z\|
\left(\sum_{|k|> n}e^{2|t|(\|V_k\| +\gamma_{bc})}\left\|
P_{k}Z^{-1}\varphi
\right\|^2\right)^{1/2},
\]
where $\gamma_{dir} = |\mathrm{Im}\,\nu|$ with $\nu$ given by \eqref{dir6} and
$\gamma_{bc} = \max\left\{|\mathrm{Im}\,\widehat p_1(0)|, |\mathrm{Im}\,\widehat p_4(0)|\right\}$ for $bc \in\{per, ap\}$.
\ec
\bpf
Follows via an application of Parseval's identity and \eqref{tildet}.
\epf

The next result  is a direct consequence of the above corollary. We cite \cite{EN00} for the standard definitions of the growth and spectral bounds of $C_0$-groups.

\bc
The growth and spectral bounds of the groups $T$ and $\widetilde T$ in Theorem \ref{mainres3} coincide.
\ec

We proceed with more explicit estimates for the groups generated by the operators $i\left(\widetilde L^P_{bc}-V\right)_{n}$, $|n|> m$. 
The case $bc = dir$ is the simplest as it yields $\left(\widetilde L^P_{bc}-V\right)_{n} = (\frac{\pi n}\omega- \nu+\  b_n)I_n$, 
where $\nu$ is given by \eqref{dir6} and 
$\{b_n: n\in\ZZ\}\in\ell^2$, according to Theorem \ref{mainres1}.

\bt\label{gdir} For $bc = dir$, the group $\widetilde T$ in \eqref{tildet} can be written as
\begeq\label{bask000}
\widetilde{T}(t)=\left(\bigoplus_{|n|>m}e^{it(\frac{\pi n}\omega -\nu + b_n)}I_n\right)\oplus e^{it\left(\widetilde L^P_{dir}-V\right)_{(m)}}
, \quad t\in\mathbb{R},
\eq
where $\{b_n: n\in\ZZ\}\in\ell^2$ can be estimated from \eqref{ltdir}.
\et

For the cases $bc \in\{per, ap\}$, we use the following formula for the group generated by a $2\times 2$ matrix \cite{K76b}:
\begeq\label{gr2}
\bs
\exp\left(it\left(
\begin{array}{cc}
a & b \\
c & d 
\end{array}
\right) 
\right) =
e^{it\frac{(a+d)}2} &\left\{
\cos(\rho t)\left(
\begin{array}{cc}
1 & 0 \\
0 & 1
\end{array}
\right)\right. \\& + \left.
\frac{i\sin(\rho t)}\rho
\left(
\begin{array}{cc}
\frac{a-d}2 & b \\
c & \frac{d-a}2 
\end{array}
\right) 
\right\},
\end{split}
\eq
where $\rho = ((a-d)^2/4+bc)^{\frac12}$, and $\frac{\sin(\rho t)}\rho$ is replaced with $t$ in the case when $\rho = 0$.

\bt\label{gperap} For $bc \in\{per, ap\}$, the group $\widetilde T$ in \eqref{tildet} can be written as
\begeq\label{bask00}
\bs
\widetilde{T}(t)&=\Bigg
(\bigoplus_{|n|>m}e^{it\left(\frac{\pi}\omega(2n+\epsilon_{bc}) -\nu + b_n^1\right)}\left\{\cos(\rho t)\left(
\begin{array}{cc}
1 & 0 \\
0 & 1
\end{array}
\right) \right. \\ & +\left.
\frac{i\sin(\rho t)}\rho
\left(
\begin{array}{cc}
\frac{\beta}2+b_n^4 & b_n^2 \\
b_n^3 & -\frac{\beta}2-b_n^4 
\end{array}
\right) 
\right\}
{\Bigg)}
\oplus e^{it\left(\widetilde L^P_{bc}-V\right)_{(m)}}
, \quad t\in\mathbb{R},
\end{split}
\eq
where $\epsilon_{bc}$ is as in \eqref{qmn}, $\beta$ -- as in \eqref{badcond}, $\rho
 = ({(\beta+2b_n^4)^2/16+b_n^2b_n^3})^{\frac12}$, $\frac{\sin(\rho t)}\rho$ is replaced with $t$ if $\rho = 0$, and
$\{b_n^k: n\in\ZZ\}\in\ell^2$, $k = 1,2,3,4$. 
\et

\bpf
Formula \eqref{bask00} is obtained via a straightforward computation. Better estimates
for the sequences $\{b_n^k: n\in\ZZ\}\in\ell^2$, $k = 1,2,3,4$, may be obtained the same way as for either Theorem \ref{mainres1000} or Theorem \ref{mainres10100}.
\epf

\brem
In our final remark, we mention a few results that may be derived from the theorems in this paper.
\ben
\item The norm of the resolvent operator $(L_{bc}-\l)^{-1}$ may be estimated.
\item For a specific potential function $P(t)$, stability of the group generated by $iL_{bc}$ may be investigated.
\item A representation of the group $\{\widetilde T(t)\}$ as a multiplicative perturbation of an explicit group may be obtained as in \cite[Corollary 2.11]{BKU18}.
\item Special cases of the results in \cite{BKR17} and \cite{BKU18} may be obtained.
\item Existence results for a semigroup in \cite{RSS99} may be improved. 
\een
\erem


\bibliographystyle{siam}
\bibliography{../refs}

\end {document}